\documentclass[a4paper,12pt]{amsart}
\newdimen\tempdimen

\def\backskip#1{%
  \settowidth{\tempdimen}{#1}%
  \multiply\tempdimen by -1
  \hskip\tempdimen
}

\def\rightalign#1{%
  \backskip{#1}%
  #1
}

\def\putleft(#1,#2)#3{%
   \put(#1,#2){\rightalign{#3}}
}

\newtheorem{theorem}{Theorem}[section]
\newtheorem{cor}[theorem]{Corollary}
\newtheorem{lemma}[theorem]{Lemma}
\newtheorem{prop}[theorem]{Proposition}
\theoremstyle{remark}

\newtheorem{remark}[theorem]{Remark}
\newtheorem{example}[theorem]{Example}

\theoremstyle{definition}
\newtheorem{definition}[theorem]{Definition}
\numberwithin{equation}{section}
\def\id{\operatorname{id}}

\def\End{\operatorname{End}}

\def\MCE{\operatorname{MCE}}
\def\cov{\operatorname{cov}}
\def\dom{\operatorname{dom}}
\def\cod{\operatorname{cod}}
\def\lsp{\operatorname{span}}
\def\clsp{\overline{\operatorname{span}}}

\newcommand{\thmref}[1]{Theorem~\ref{#1}}

\newcommand{\lemref}[1]{Lemma~\ref{#1}}
\newcommand{\propref}[1]{Prop\-o\-si\-tion~\ref{#1}}
\newcommand{\defref}[1]{Def\-i\-ni\-tion~\ref{#1}}
\newcommand{\corref}[1]{Cor\-ol\-lary~\ref{#1}}
\newcommand{\field}[1]{\mathbb{#1}}
\newcommand{\CC}{\field{C}}

\newcommand{\NN}{\field{N}}
\newcommand{\RR}{\field{R}}

\newcommand{\ZZ}{\field{Z}}

\newcommand{\Bb}{{B}}

\newcommand{\Hh}{{\mathcal H}}

\newcommand{\Kk}{{\mathcal K}}
\newcommand{\Ll}{{\mathcal L}}

\newcommand{\Oo}{{\mathcal O}}

\newcommand{\Tt}{{\mathcal T}}

\begin{document}
\title[Product systems of graphs]{\boldmath Product systems of graphs and\\
the Toeplitz algebras of higher-rank graphs}
\author[Iain Raeburn]{Iain Raeburn}
\author[Aidan Sims]{Aidan Sims}
\address{School of Mathematical and Physical Sciences\\
University of Newcastle\\
NSW 2308\\
Australia}
\email{iain, aidan@frey.newcastle.edu.au}
\date{March 21, 2002}

\thanks{This research was supported by the Australian Research Council.}

\begin{abstract}
There has recently been much interest in the $C^*$-algebras of
directed graphs. Here
we consider product systems $E$ of directed graphs over
semigroups and associated
$C^*$-algebras $C^*(E)$ and $\mathcal{T}C^*(E)$ which generalise the
higher-rank graph
algebras of Kumjian-Pask and their Toeplitz analogues.
  We study these algebras by constructing from $E$ a product system
$X(E)$ of Hilbert
bimodules, and applying recent  results of Fowler about the
Toeplitz algebras of such systems. Fowler's hypotheses turn out to be
very interesting
graph-theoretically, and indicate new relations which will have to be
added to the
usual Cuntz-Krieger relations to obtain a satisfactory theory of Cuntz-Krieger
algebras for product systems of graphs; our
algebras
$C^*(E)$ and $\mathcal{T}C^*(E)$ are universal for families of
partial isometries
satisfying these relations.

Our main result is a uniqueness theorem for $\mathcal{T}C^*(E)$ which
has particularly
interesting implications for the $C^*$-algebras of non-row-finite
higher-rank graphs.
  This theorem is apparently beyond the reach of
Fowler's theory, and  our proof requires a detailed analysis of the
expectation onto the diagonal in $\mathcal{T}C^*(E)$.
\end{abstract}

\maketitle

\section{Introduction}\label{intro}

The $C^*$-algebras $C^*(E)$ of infinite directed graphs $E$ are
generalisations of
the Cuntz-Krieger algebras which include many interesting
$C^*$-algebras and provide a
rich supply of models for simple purely infinite algebras (see, for example,
\cite{KPRR, D, HS, Sz}). In the first papers, it was assumed for
technical reasons
that the graphs were locally finite. However, after $C^*(E)$ had been
realised as the
Cuntz-Pimsner algebra $\Oo_{X(E)}$ of a Hilbert bimodule $X(E)$ in
\cite{FR}, it was noticed that $\Oo_{X(E)}$ made sense for arbitrary
infinite graphs.
The analysis in \cite{FR} applied to the Toeplitz algebra $\Tt_{X(E)}$ rather
than $\Oo_{X(E)}$, but the two coincide for some infinite graphs $E$, and
hence the results of \cite{FR} gave information about $\Oo_{X(E)}$ for these
graphs. The results of \cite{FR} therefore suggested an appropriate
definition of
$C^*(E)$  for arbitrary
$E$, which was implemented in
\cite{FLR}.

Higher-rank analogues of Cuntz-Krieger algebras  and of the $C^*$-algebras of
row-finite graphs have been studied by Robertson-Steger \cite{RS} and
Kumjian-Pask
\cite{KP}, respectively. It was observed in \cite{FS} that the
higher-rank graphs of
Kumjian and Pask could be viewed as product systems of graphs over
the semigroup
$\NN^k$. The main object of this paper is to extend the construction
$E\mapsto X(E)$ to product systems of graphs over $\NN^k$ and other
semigroups, to
apply the results of
\cite{F99} to the resulting product systems of Hilbert bimodules, and
to see what
insight might be gained into the $C^*$-algebras of arbitrary
higher-rank graphs.

It is relatively easy to extend the construction of $X(E)$ to product
systems, and to identify  \emph{Toeplitz
$E$-families} which correspond to the Toeplitz representations of $X(E)$
studied in \cite{F99}. The story becomes interesting when we
investigate the conditions on $E$ and on Toeplitz
$E$-families which ensure that we can apply \cite[Theorem~7.2]{F99} to the
corresponding representation of $X(E)$. To understand the issues, we digress
briefly.

The isometric representation theory of semigroups suggests that in general
$\Tt_{X(E)}$ will be too big to behave like a Cuntz-Krieger algebra,
and that we
should restrict attention to the
Nica-covariant representations of \cite{N92, LR, F, F99}. However, Nica
covariance is in general a spatial phenomenon, and to talk about the universal
$C^*$-algebra $\Tt_{\cov}(X)$ generated by a Nica-covariant Toeplitz
representation
of a product system $X$ of bimodules, we need to assume that
$X$ is compactly aligned in the sense of \cite{F, F99}.

We identify the \emph{finitely aligned}
product systems $E$ of graphs for which $X(E)$ is compactly aligned, and the
\emph{Toeplitz-Cuntz-Krieger $E$-families}
$\{S_\lambda\}$ which correspond to Nica-covariant Toeplitz
representations of $X(E)$.
The $C^*$-algebra generated by
$\{S_\lambda\}$ is then spanned by the products $S_\lambda S_\mu^*$,  as
Cuntz-Krieger algebras and their Toeplitz analogues are. We therefore
define the Toeplitz algebra $\Tt C^*(E)$ of a finitely aligned product system
$E$ to be the universal $C^*$-algebra generated by a Toeplitz-Cuntz-Krieger
$E$-family; for technical reasons, we only define the
Cuntz-Krieger algebra
$C^*(E)$ to be the appropriate quotient of
$\Tt C^*(E)$ when $E$ has no sinks.

Fowler's \cite[Theorem~7.2]{F99}
gives a spatial condition under which a Nica-covariant Toeplitz
representation of
a compactly aligned  product system
$X$ of Hilbert bimodules is faithful on
$\Tt_{\cov}(X)$. Since $\Tt C^*(E)$ has essentially the same
representation theory as
$\Tt_{\cov}(X(E))$, Fowler's theorem describes some faithful
representations of $\Tt C^*(E)$. However, the resulting theorem about
Toeplitz-Cuntz-Krieger
$E$-families is not as sharp as we would like, for the same reasons that
\cite[Theorem~2.1]{FR} is not: applied to the single graph
$E$ with $\Tt C^*(E)=\Oo_\infty$, it says that isometries $\{S_i\}$
satisfying $1>\sum_{i=1}^\infty S_iS_i^*$ generate an isomorphic copy
of $\Oo_\infty$, whereas we know from  \cite{C} that $1\geq \sum_{i=1}^\infty
S_iS_i^*$ suffices. Our main theorem
is sharp in this sense: it is an analogue of
\cite[Theorem~3.1]{FR} rather than \cite[Theorem~2.1]{FR}. It suggests
an appropriate set of Cuntz-Krieger relations for product systems of
not-necessarily-row-finite graphs, and
gives a uniqueness theorem of Cuntz-Krieger type for
$k$-graphs in which each vertex receives infinitely many edges of each
degree.

  We start with a short
review of the basic facts about graphs and  the
Cuntz-Krieger bimodule
$X(E)$ of a single graph $E$. In \S\ref{sec1}, we associate to each
product system
$E$ of graphs  a product system $X(E)$ of
Cuntz-Krieger bimodules (Proposition~\ref{prop:{X(E_p)} a PS}). In
\S\ref{sec2}, we define Toeplitz $E$-families, and show that there is a
one-to-one correspondence between such families and Toeplitz representations of
$X(E)$ (Theorem~\ref{thm:families and repns}). We then restrict attention to
product systems over the quasi-lattice ordered semigroups of Nica,
and identify the
finitely aligned product systems $E$ of graphs for which
$X(E)$ is compactly aligned (Theorem~\ref{thm:compactly aligned}). In
\S\ref{sec4}, we
discuss Nica covariance, and show that for finitely aligned systems,
it becomes a familiar relation which is automatically satisfied by
Cuntz-Krieger
families of a single graph. By adding this relation to those of a
Toeplitz family, we
obtain an appropriate definition of Toeplitz-Cuntz-Krieger
$E$-families for more general $E$, and then
$\Tt C^*(E)$ is the universal $C^*$-algebra generated by such a
family. We can now
apply Fowler's theorem to $X(E)$
(Proposition~\ref{thm:7.2}), and
deduce that the Fock representation of $\Tt C^*(E)$ is faithful
(Corollary~\ref{cor:Fock faithful}).

Our main Theorem~\ref{thm:faithfulness} is a $C^*$-algebraic
uniqueness theorem.
It does not appear to follow from Fowler's results: its proof
requires a detailed
analysis of the expectation onto the diagonal in $\Tt C^*(E)$ and its spatial
implementation, as well as an application of Corollary~\ref{cor:Fock
faithful}.  In
the last section,
  we apply Theorem~\ref{thm:faithfulness} to the $k$-graphs of \cite{KP}. Our
results are all interesting in this case, and those interested primarily in
$k$-graphs could assume $P=\NN^k$ throughout the paper without losing
the main points.

\section{Preliminaries}\label{prelims}

\subsection{Graphs and Cuntz-Krieger families}\label{sec:P1}
A {\em directed graph} $E = (E^0, E^1, r, s)$ consists of a countable
vertex set $E^0$, a countable edge set $E^1$, and range and source
maps $r, s:E^1\to E^0$. All graphs in this paper are
directed.

A {\em Toeplitz-Cuntz-Krieger $E$-family} in a
$C^*$-algebra $B$ consists  of
mutually orthogonal projections $\{p_v : v \in E^0\}$ in $B$ and
partial isometries $\{s_\lambda : \lambda \in E^1\}$ in $B$
satisfying $s^*_{\lambda}s_{\lambda} = p_{r(\lambda)}$ for
$\lambda \in E^1$ and
\[
p_v \ge \sum_{\lambda \in F} s_{\lambda}s^*_{\lambda} \quad\text{for
every $v \in E^0$ and every finite set $F \subset s^{-1}(v)$}.
\]
It is a {\em Cuntz-Krieger $E$-family} if
\[
p_v = \sum_{\lambda \in s^{-1}(v)}
s_{\lambda}s^*_{\lambda} \quad\text{whenever $s^{-1}(v)$ is finite
and nonempty}.
\]

\subsection{Hilbert bimodules}\label{sec:P2}

Let $A$ be a $C^*$-algebra. A {\em right-Hilbert $A-A$ bimodule} (or {\em
Hilbert bimodule over $A$}) is a right Hilbert $A$-module $X$ together
with a left action $(a,x)\mapsto
a\cdot x$ of
$A$ by adjointable operators on $X$; we denote by $\phi$ the homomorphism of
$A$ into $\Ll(X)$ given by the left action. We say
$X$ is {\em essential} if
\[
\clsp\{a \cdot x : a \in A,\; x \in X\} = X.
\]

A {\em Toeplitz
representation} $(\psi,\pi)$ of a Hilbert bimodule $X$ in a $C^*$-algebra $B$
consists of a linear map $\psi : X \to B$ and a homomorphism
$\pi: A \to B$ such that
\[
\psi(x\cdot a) = \psi(x)\pi(a),\quad
\psi(a \cdot x) = \pi(a)\psi(x),\quad\text{and}\quad
\psi(x)^*\psi(y) = \pi(\langle x,y
\rangle_A)
\]
for $x,y \in X$ and $a \in A$. There is then a unique homomorphism
$\psi^{(1)} : \Kk(X)
\to B$ such that
\[
\psi^{(1)}(\Theta_{x,y}) = \psi(x)\psi(y)^* \quad\text{for $x,y \in
X$;}
\]
see \cite[page~202]{P97}, \cite[Lemma~2.2]{KPW}, or
\cite[Remark~1.7]{FR} for details. The representation
$(\psi, \pi)$ is {\em Cuntz-Pimsner covariant} if
\[
\psi^{(1)}(\phi(a)) = \pi(a) \quad\text{whenever $\phi(a) \in
\Kk(X)$.}
\]

Pimsner associated to each Hilbert bimodule $X$
a $C^*$-algebra $\Tt_X$ which is universal
for Toeplitz representations of $X$, and a quotient $\Oo_X$ which
is universal for Cuntz-Pimsner covariant Toeplitz representations of $X$
(\cite{P97}; see also \cite[\S1]{FR}).

\subsection{Cuntz-Krieger bimodules}
The Cuntz-Krieger bimodule $X(E)$ of a graph $E$, as in
\cite[Example~1.2]{FR}, consists of the functions
$x: E^1
\to
\CC$ such that
\begin{equation}\label{eq:rho-x}
\rho_x:v\mapsto \sum_{\lambda \in E^1,
r(\lambda)=v} \lvert x(\lambda) \rvert^2
\end{equation}
vanishes at infinity on $E^0$. With
\begin{gather*}
(x\cdot a)(\lambda) := x(\lambda)a(r(\lambda)) \text{ and }
(a\cdot x)(\lambda) := a(s(\lambda))x(\lambda) \quad\text{for
$\lambda \in E^1$, and}\\
\langle x,y \rangle_{C_0(E^0)}(v) := \sum_{\lambda \in E^1,
r(\lambda)=v} \overline{x(\lambda)}y(\lambda) \quad\text{for $v \in
E^0$,}
\end{gather*}
$X(E)$ is a Hilbert bimodule over
$C_0(E^0)$. The Toeplitz representations of $X(E)$ are in one-to-one
correspondence
with the Toeplitz-Cuntz-Krieger $E$-families via
$(\psi,\pi) \leftrightarrow \{\psi(\delta_\lambda),\pi(\delta_v)\}$
\cite[Example~1.2]{FR}. Hence $\Tt_{X(E)}$ is universal for
Toeplitz-Cuntz-Krieger
$E$-families. When $E$ has no sinks, the left action of $C_0(E^0)$ on
$X(E)$ is faithful, the Cuntz-Pimsner covariant representations correspond to
Cuntz-Krieger
$E$-families, and the quotient
$\Oo_{X(E)}$ is the usual graph $C^*$-algebra $C^*(E)$.

Because of the correspondence $(\psi,\pi) \leftrightarrow
\{\psi(\delta_\lambda),\pi(\delta_v)\}$, it is convenient in
calculations to work with the point masses $\delta_\lambda\in X(E)$.
The following
lemma explains why this suffices.

\begin{lemma}\label{lem:xcdense}
The space $X_c(E):=C_c(E^1)$ is a dense submodule of $X(E)$, and the point
masses $\{\delta_\lambda:\lambda\in E^1\}$ are a vector-space basis for
$X_c(E^1)$.
\end{lemma}

\begin{proof}
As a Banach space, $X(E)$ is the $c_0$-direct sum $\bigoplus_{v\in
E^0}\ell^2(r^{-1}(v))$, and $X_c(E)$ is the algebraic direct sum of
the subspaces $C_c(r^{-1}(v))$. So it is standard that
$X_c(E)$ is dense. For $x\in X_c(E)$, we have $x=\sum_{\lambda\in
E^1}x(\lambda)\delta_\lambda$.
\end{proof}

\section{Product systems of graphs and of Hilbert
bimodules}\label{sec1}
Throughout the next two sections, $P$ denotes an arbitrary countable
semigroup with identity $e$. If $E = (E^0, E^1, r_E, s_E)$ and $F = (E^0, F^1,
r_F, s_F)$ are two graphs with the same vertex set $E^0$,  then $E
\times_{E^0} F$ denotes the graph with $(E \times_{E^0} F)^0 :=
E^0$,
\[
(E \times_{E^0} F)^1 := \{(\lambda,\mu) : \lambda \in E^1, \mu \in
F^1, r_E(\lambda) = s_F(\mu)\},
\]
and $s(\lambda,\mu) := s_E(\lambda)$, $r(\lambda,\mu) := r_F(\mu)$.

We recall from \cite{FS} that a {\em product system $(E,\varphi)$ of graphs
over $P$} consists of graphs $\{(E^0, E^1_p,
r_p, s_p) : p \in P\}$ with common vertex set $E^0$ and disjoint edge sets
$E_p^1$, and  isomorphisms
$\varphi_{p,q} : E_p \times_{E^0} E_q \to E_{pq}$ for $p, q\in P$
satisfying the associativity condition
\begin{equation}\label{eqn:ps associative}
\varphi_{pq,r}(\varphi_{p,q}(\lambda,\mu)),\nu)
= \varphi_{p,qr}(\lambda,\varphi_{q,r}(\mu,\nu))
\end{equation}
for all $p,q,r \in P$, $(\lambda,\mu) \in (E_p \times_{E^0} E_q)^1$, and
$(\mu,\nu)\in (E_q \times_{E^0} E_r)^1$; we require that
\[
E_e = (E^0,E^0,\id_{E^0},\id_{E^0}).
\]
We write $d(\lambda)=p$ to mean $\lambda\in E^1_p$;
because the $E_p^1$ are disjoint, this gives a well-defined \emph{degree map}
$d:E^1:=\bigcup_{p\in P} E_p^1\to P$, which gives the vertices
$E^0=E_e^1$ degree
$e$. The range and source maps combine to give maps $r,s:E^1\to E^0$.

The isomorphisms $\varphi_{p,q}$ in a product system $(E,\varphi)$ combine
to give a partial multiplication on $E^1$: for $(\lambda,\mu)\in
E^1_p\times_{E^0}E^1_q$, we define
$\lambda\mu=\varphi_{p,q}(\lambda,\mu)\in E^1_{pq}$. This multiplication
is associative by \eqref{eqn:ps associative}.  Since each $\varphi_{p,q}$
is an isomorphism, the multiplication has the following \emph{factorisation
property}: for each $\gamma \in E^1_{pq}$, there is a unique
$(\lambda,\mu) \in (E_p\times_{E^0} E_q)^1$ such that $\gamma =
\lambda\mu$. It follows that if $\lambda \in E^1_{pqr}$, then there is a
unique $\lambda(p, pq)\in E_q^1$ such that  $\lambda=\lambda' \lambda(p, pq)
\lambda''$ with $d(\lambda') = p$ and
$d(\lambda'') = r$. By \eqref{eqn:ps associative} and the factorisation
property, $s(\lambda)\lambda=\lambda=\lambda r(\lambda)$ for all $\lambda$.

A single graph $E$ gives a product system over $\NN$ in which $E_n^1$
consists of the
paths of length $n$ in $E$. More generally:

\begin{example}[$k$-graphs]\label{k-graphs and PSs}
It is shown in \cite[Examples~1.5, (4)]{FS} that the product
systems of graphs over $\NN^k$ are essentially the same as the $k$-graphs of
\cite[Definitions~1.1]{KP}:
\begin{itemize}
\item[$\bullet$]
Given a product system $(E, \varphi)$ of graphs over $\NN^k$, let
$\Lambda_E$ be the category with objects $E^0$ and morphisms $E^1$, with
$\dom(\lambda) := r(\lambda)$ and
$\cod(\lambda) := s(\lambda)$. The degree map is that of $E$, the
morphism $\lambda\circ\mu$ is by definition the morphism associated to the
edge $\lambda\mu$, and the factorisation property for
$\Lambda_E$ reduces to that of $E$.

\item[$\bullet$]
Given a $k$-graph $(\Lambda, d)$, let
$(E_\Lambda)^0 := \Lambda^0$, $(E_\Lambda)^1_n := \Lambda^n$ for
$n \in \NN^k$, $\lambda\mu := \lambda\circ\mu \in \Lambda^{m+n}$
whenever $(\lambda,\mu) \in (E_m \times_{E^0} E_n)^1$, and define
$r:=\dom$ and $s:=\cod$.
\end{itemize}
The direction of the edges is reversed in going from $(\Lambda, d)$ to
$(E_\Lambda, \varphi_\Lambda)$ to ensure that the representations of
the two coincide (compare \defref{def:TCK family} with
\cite[Definitions~1.5]{KP}).
\end{example}

\begin{prop} \label{prop:{X(E_p)} a PS}
If $(E, \varphi)$ is a product system of graphs over $P$, then there
is a unique associative multiplication on
$X(E):=\bigcup_{p\in P}X(E_p)$ such that
\begin{equation}\label{eq:multiplication}
\delta_\lambda \delta_\mu :=
\begin{cases}
\delta_{\lambda\mu} &\text{if $(\lambda,\mu) \in (E_{d(\lambda)}
\times_{E^0} E_{d(\mu)})^1$} \\
0 &\text{otherwise,}
\end{cases}
\end{equation}
and $X(E)$ thus becomes a product system of Hilbert
bimodules over $C_0(E^0)$ as in \cite[Definition~2.1]{F99}.
\end{prop}

\begin{remark}
We have described the multiplication using point masses because
we want to use them in
calculations. However, we also write it out explicitly in
Corollary~\ref{cor:explicitm}.
\end{remark}

\begin{proof}[Proof of Proposition~\ref{prop:{X(E_p)} a PS}]
It follows from Lemma~\ref{lem:xcdense} that the elements
$\delta_\lambda\otimes\delta_\mu$ are a basis for the algebraic tensor
product $X_c(E_p)\odot X_c(E_q)$, and hence there is a well-defined
linear map $\pi:X_c(E_p)\odot X_c(E_q)\to X_c(E_{pq})$ such that
\[
\pi(\delta_\lambda \otimes\delta_\mu) =
\begin{cases}
\delta_{\lambda\mu} &\text{if $(\lambda,\mu) \in (E_{d(\lambda)}
\times_{E^0} E_{d(\mu)})^1$} \\
0 &\text{otherwise.}
\end{cases}
\]
Let $\lambda,\mu,\eta,\xi\in E^1$. Then
\begin{align}
\big\langle \delta_\lambda \otimes \delta_\mu, \delta_\eta \otimes
\delta_\xi \big\rangle_{C_0(E^0)}(v)
&= \big\langle \langle\delta_\eta, \delta_\lambda\rangle_{C_0(E^0)} \cdot
\delta_\mu, \delta_\xi \big\rangle_{C_0(E^0)}(v) \notag\\
&=
\begin{cases}
1 &\text{if $\eta = \lambda$, $\xi = \mu$, $r(\lambda) = s(\mu)$ and
$r(\mu) = v$} \\
0 &\text{otherwise.}
\end{cases}\label{2.5}
\end{align}
On the other hand,
\begin{align*}
\big\langle \pi(\delta_\lambda\otimes\delta_\mu) ,
\pi(\delta_\eta\otimes\delta_\xi)
&\big\rangle_{C_0(E^0)}(v)\\
&=
\begin{cases}
\langle \delta_{\lambda\mu} , \delta_{\eta\xi} \rangle_{C_0(E^0)}(v)
&\text{if $r(\lambda) = s(\mu)$ and $r(\eta) = s(\xi)$} \\
0 &\text{otherwise}
\end{cases} \\
&=
\begin{cases}
1 &\text{if $r(\lambda) = s(\mu)$, $r(\eta) = s(\xi)$, $\lambda\mu
= \eta\xi$ and $r(\mu)=v$}\\
0 &\text{otherwise,}
\end{cases}
\end{align*}
which by the factorisation property is (\ref{2.5}).
Since $X_c(E_p)$ is dense in $X(E_p)$ (Lemma~\ref{lem:xcdense}), it
follows that $\pi$
extends to an isometric linear isomorphism of
$X(E_p)\otimes_{C_0(E^0)}X(E_q)$ onto
$X(E_{pq})$. It is easy to check on dense subspaces $X_c(E_p)$ and
$\lsp\{\delta_v\}\subset C_0(E^0)$ that
$\pi$ is an isomorphism of Hilbert $C_0(E^0)$-bimodules. We
now define $xy:=\pi(x\otimes y)$, and associativity of
this multiplication follows from \eqref{eqn:ps
associative}. More calculations on dense subspaces show that $xa=x\cdot
a$ and $ax=a\cdot x$ for $a\in C_0(E^0)=X(E_e)$ and $x\in X(E_p)$.
\end{proof}

\begin{cor}\label{cor:explicitm}
For $x\in X(E_p)$ and $y\in X(E_q)$, we have
\begin{equation}\label{formforprod}
(xy)(\lambda\mu)=x(\lambda)y(\mu)\quad\text{for $(\lambda,\mu)\in
(E_p\times_{E^0}E_q)^1$.}
\end{equation}
\end{cor}

\begin{proof}
The multiplication extends to an isomorphism of
$X(E_p)\otimes_{C_0(E^0)}X(E_q)$ onto $X(E_{pq})$, $(x,y)\mapsto
x\otimes y$ is continuous, and the various evaluation maps $z\mapsto
z(\lambda)$ are continuous, so Lemma~\ref{lem:xcdense} implies that it
is enough to prove (\ref{formforprod}) for
$x\in X_c(E_p)$ and $y\in X_c(E_q)$. For such $x,y$ we have
\[
(xy)(\lambda\mu)=\sum_{\alpha\in E_p^1,\beta\in E_q^1}
x(\alpha)y(\beta)(\delta_\alpha\delta_\beta)(\lambda\mu),
\]
which collapses to $x(\lambda)y(\mu)$ by the factorisation property.
\end{proof}

\section{Representations of product systems}\label{sec2}
Throughout this section, $(E,\varphi)$ is a product system of graphs
over $P$.

\begin{definition}\label{def:TCK family}
Partial isometries $\{s_\lambda : \lambda \in E^1\}$ in a
$C^*$-algebra $B$ form a
{\em Toeplitz $E$-family} if:
\begin{itemize}
\item[(1)] $\{s_v : v \in E^0\}$ are mutually orthogonal
projections,
\item[(2)] $s_\lambda s_\mu = s_{\lambda\mu}$ for all $\lambda,\mu
\in E^1$ such that $r(\lambda) = s(\mu)$,
\item[(3)] $s^*_\lambda s_\lambda = s_{r(\lambda)}$ for all $\lambda
\in E^1$, and
\item[(4)] for all $p \in P\setminus\{e\}, v \in E^0$ and every
finite $F \subset s^{-1}_p(v)$, $s_v \ge \sum_{\lambda
\in F} s_\lambda s^*_\lambda$.
\end{itemize}
\end{definition}

We recall from \cite{F99} that a Toeplitz representation $\psi$ of a
product system
$X$ of bimodules consists of linear maps $\psi_p:X_p\to B$ such that each
$(\psi_p,\psi_e)$ is a Toeplitz representation of $X_p$, and
$\psi_p(x)\psi_q(y)=\psi_{pq}(xy)$. It is Cuntz-Pimsner covariant if each
$(\psi_p,\psi_e)$ is Cuntz-Pimsner covariant. Fowler proves that there is a
$C^*$-algebra $\Tt_X$ generated by a universal Toeplitz
representation $i_X$, and a
quotient $\Oo_X$ generated by a universal Cuntz-Pimsner covariant
representation $j_X$
\cite[\S2]{F99}.

\begin{theorem}\label{thm:families and repns}
Let $(E, \varphi)$ be a product system of graphs over a semigroup $P$, and let
$X(E)$ be the corresponding product system of Cuntz-Krieger bimodules. If
$\psi$ is a Toeplitz representation of $X(E)$, then
\begin{equation}\label{family from repn}
\{s_\lambda := \psi_{d(\lambda)} (\delta_\lambda):\lambda \in E^1\}
\end{equation}
is a Toeplitz $E$-family; conversely, if
$\{s_\lambda : \lambda \in E^1\}$ is a Toeplitz
$E$-family, then the map
\begin{equation}\label{repn from family}
x \in C_c(E^1_p) \mapsto \sum_{\lambda \in E^1_p} x(\lambda)s_\lambda
\end{equation}
extends to a Toeplitz representation of $X(E)$ from which we can recover
$s_\lambda =
\psi_{d(\lambda)} (\delta_\lambda)$. The representation $\psi$ is Cuntz-Pimsner
covariant if and only if $\{s_\lambda\}$ satisfies
\begin{equation} \label{eq:CP rel}
s_v = \sum_{\lambda \in
s^{-1}_p(v)} s_\lambda s^*_\lambda \text{ whenever $s^{-1}_p(v)$ is
finite (possibly
empty).}
\end{equation}
\end{theorem}

\begin{proof}
If $\psi$ is a Toeplitz representation of $X(E)$, then
\cite[Example~1.2]{FR} shows that
\[
\{\psi_e(\delta_v),
\psi_p(\delta_\lambda) : v \in E^0, \lambda \in E^1_p\}
\]
is a Toeplitz-Cuntz-Krieger family for
$E_p$ as in \cite{FR}, and this gives (1), (3), and (4) of \defref{def:TCK
family}. \defref{def:TCK family}(2) follows from (\ref{eq:multiplication})
because $\psi$ is a homomorphism.

Now suppose that $\psi$ is Cuntz-Pimsner covariant and $s_p^{-1}(v)$
is finite. Write
$\phi_p : C_0(E^0) \to \Ll(X_p)$ for the homomorphism that implements the left
action on $X_p$. Then
\begin{equation}\label{eqn:CP-cov}
\sum_{\lambda \in s^{-1}_p(v)} \psi_p(\delta_\lambda)
\psi_p(\delta_\lambda)^*
= \sum_{\lambda \in s^{-1}_p(v)} \psi^{(1)}_p(\Theta_{\delta_\lambda,
\delta_\lambda})
= \psi^{(1)}_p \Big( \sum_{\lambda \in s^{-1}_p(v)}
\Theta_{\delta_\lambda, \delta_\lambda}\Big).
\end{equation}
For $x \in X_p$, $w \in E^0$ and $\mu \in E^1_p$,
\[
\Big(\sum_{\lambda \in s^{-1}_p(w)} \Theta_{\delta_\lambda,
\delta_\lambda}(x)\Big)(\mu)
= \left\{\hskip-5pt\begin{array}{ll}
x(\mu) &\text{if $\mu \in s^{-1}_p(w)$} \\
0 &\text{otherwise}
\end{array}\hskip-5pt\right\}
= (\delta_w \cdot x)(\mu).
\]
Hence the right hand side of \eqref{eqn:CP-cov} is just
$\psi^{(1)}_p(\phi_p(\delta_v))$. Since
$\phi_p(\delta_v)$ belongs to $\Kk(X_p)$ \cite[Proposition~4.4]{FR},
Cuntz-Pimsner
covariance gives
$\psi^{(1)}_p(\phi_p(\delta_v)) = \psi_e(\delta_v)$. Thus
\[
\sum_{\lambda \in s^{-1}_p(v)} s_\lambda s_\lambda^* = \sum_{\lambda
\in s^{-1}_p(v)} \psi_p(\delta_\lambda) \psi_p(\delta_\lambda)^* =
\psi_e(\delta_v) = s_v.
\]

If $\{s_\lambda : \lambda \in E^1\}$ is a
Toeplitz $E$-family, \cite[Example~1.2]{FR} implies
that $\psi_p(\delta_\lambda) := s_\lambda$ extend to Toeplitz
representations $(\psi_p, \psi_e)$ of $X_p$ for $p \in P$; since
\[
\psi_{pq}(\delta_{\lambda}\delta_\mu)=\psi_{pq}(\delta_{\lambda\mu})=s_{\lambda\
mu}=s_\lambda
s_\mu=\psi_p(\delta_\lambda)\psi_q(\delta_\mu),
\]
it follows that
$\psi$
is a Toeplitz representation of $X(E)$. We trivially have $s_\lambda =
\psi_{d(\lambda)} (\delta_\lambda)$.

If $\{s_\lambda : \lambda \in E^1\}$ satisfies (\ref{eq:CP rel}),
then for $p \in P$ and $v \in E^0$ with $s^{-1}_p(v)$ finite,
\[
\psi^{(1)}_p(\phi_p(\delta_v))
= \psi^{(1)}_p \Big(\sum_{\lambda\in s^{-1}_p(v)}
\Theta_{\delta_\lambda, \delta_\lambda}\Big)
= \sum_{\lambda\in s^{-1}_p(v)} \psi_p(\delta_\lambda)
\psi_p(\delta_\lambda)^*,
\]
which is $\psi_e(\delta_v)$ by  (\ref{eq:CP rel}).
Proposition~4.4 of \cite{FR} ensures that $\{\delta_v : \lvert
s^{-1}_p(v) \rvert
< \infty\}$ spans a dense subspace of $\{a\in C_0(E^0):\phi(a)\in
\Kk(X_p)\}$, so
$\psi$ is Cuntz-Pimsner covariant.
\end{proof}

\begin{cor} \label{OoX and TtX univ for graph}
Let $(E,\varphi)$ be a product system of graphs over a semigroup $P$. Then
$(\Tt_{X(E)}, i_{X(E)})$ is universal for Toeplitz
$E$-families in the sense that
\begin{itemize}
\item[(1)] $\{s_\lambda\} := \{i_{X(E)}(\delta_\lambda)\}$ is a
Toeplitz $E$-family which generates $\Tt_{X(E)}$; and
\item[(2)] for every Toeplitz $E$-family $\{s_\lambda\}$, there is a
representation
$\psi_*$ of $\Tt_{X(E)}$ such that $(\psi_* \circ i_{X(E)})(\delta_\lambda) =
s_\lambda$ for every
$\lambda \in E^1$.
\end{itemize}
Similarly, $(\Oo_{X(E)}, j_{X(E)})$ is universal for
Toeplitz $E$-families satisfying \eqref{eq:CP rel}.
\end{cor}

\begin{proof}
This follows from Theorem~\ref{thm:families and repns} and the
universal properties
of $\Tt_{X(E)}$ and $\Oo_{X(E)}$ described in \cite[Propositions~2.8
and 2.9]{F99}.
\end{proof}

If $(E, \varphi)$ is a product system of row-finite graphs without sinks over
$\NN^k$, then $\Lambda_E$ is row-finite and has no sources
as in \cite{KP}, and the Toeplitz
$E$-families which satisfy \eqref{eq:CP rel} are precisely the
$*$-representations of
$\Lambda_E$. Hence:

\begin{cor} \label{cor:C*(Lambda) = OoX}
Let $\Lambda$ be a row-finite $k$-graph with no sources as in
\cite{KP}, define $E_\Lambda$  as in Example~\ref{k-graphs and PSs}, and let
$X = X(E_\Lambda)$. Then there is an isomorphism of $C^*(\Lambda)$ onto $\Oo_X$
carrying
$s_\lambda$ to $i_X(\delta_\lambda)$.
\end{cor}

\begin{remark}\label{bloodysinks}
If there are vertices which are sinks in one or more $E_p$, then some
subtle issues arise, and the Toeplitz $E$-families satisfying
\eqref{eq:CP rel} are
not necessarily the Cuntz-Krieger $\Lambda_E$-families studied in
\cite{RSY}.  Here, though, we care primarily about Toeplitz
familes, and the presence of sinks does not cause problems.
\end{remark}

\section{Compactly aligned product systems of Cuntz-Krieger
bimodules}\label{sec3}

The compactly aligned product systems are a large class of product
systems whose Toeplitz algebras have
been analysed in \cite{F} and \cite{F99}. To apply the results of \cite{F99},
we need to identify the product systems $E$ of graphs for which $X(E)$ is
compactly aligned.

In compactly aligned product systems, the underlying semigroup $P$ has to be
quasi-lattice ordered in the sense of Nica
\cite{N92, LR}. Suppose $P$ is a subsemigroup of a group $G$ such
that $P\cap P^{-1} = \{e\}$. Then $g \le h \Longleftrightarrow
g^{-1}h \in P$ defines a partial order on $G$, and $P$ is \emph{quasi-lattice
ordered} if every finite subset of $G$ with an upper bound in $P$ has
a least upper
bound in $P$. (Strictly speaking, it is the pair $(G,P)$ which is quasi-lattice
ordered.) If two
elements $p$ and $q$ have a
common upper bound in $P$, $p \vee q$ denotes their least upper bound;
otherwise, we write $p\vee q = \infty$.

Totally ordered groups, free groups, and products of these groups are all
quasi-lattice ordered. The main example of interest to us is
$(G,P)=(\ZZ^k,\NN^k)$, which is actually \emph{lattice-ordered}: each pair
$m,n \in \NN^k$ has a least upper bound $m\vee n$ with $i$th coordinate $(m
\vee n)_i :=\max\{m_i,n_i\}$.

Let $X$ be a product system of bimodules over a quasi-lattice ordered
semigroup $P$, and suppose $p,q \in P$ have $p\vee q < \infty$. Since
$S\in \Ll(X_p)$ acts as an adjointable operator $S\otimes 1$ on $X_p
\otimes_A X_{p^{-1}(p \vee q)}$, the isomorphism
of $X_p
\otimes_A X_{p^{-1}(p \vee q)}$ onto $X_{p \vee q}$ induced by the
multiplication gives an action of $\Ll(X_p)$ on $X_{p \vee q}$; we
write $S^{p\vee q}_p$ for the image of $S\in\Ll(X_p)$, so that $S^{p\vee
q}_p$ is characterised by
\begin{equation}\label{defSpq}
S^{p\vee q}_p(xy) := (Sx)y\quad\text{ for $x\in X_p$, $y\in X_{p^{-1}(p
\vee q)}$.}
\end{equation}
The product system $X$ is \emph{compactly aligned}
\cite[Definition~5.7]{F99} if
\[
S \in \Kk(X_p)\text{ and }T\in\Kk(X_q)\text{ imply } (S^{p\vee
q}_p)(T^{p\vee q}_q) \in \Kk(X_{p \vee q}).
\]
When $X = X(E)$ is a product system of
Cuntz-Krieger  bimodules, Lemma~\ref{lem:xcdense} implies that the
point masses span
dense subspaces of
$X(E_p)$, and the rank-one operators $\Theta_{x,y}$ span dense subspaces of
$\Kk(X)$; thus to prove that
$X(E)$ is compactly aligned, it suffices to check that every
\begin{equation}\label{eq:rank one}
(\Theta_{\delta_{\mu_1}, \delta_{\mu_2}})^{p\vee
q}_p(\Theta_{\delta_{\nu_1}, \delta_{\nu_2}})^{p\vee q}_q \quad\text{belongs
to $\Kk(X(E_{p\vee q}))$.}
\end{equation}
To prove that a given $X(E)$ is not compactly aligned, we need
to be able to recognise non-compact operators on $X(E)$.

\begin{lemma}\label{lem:compacts C_0}
Let $X(E)$ be the Cuntz-Krieger bimodule of a graph, and let $S \in
\Kk(X(E))$. Then the function $x_S:E^1\to \RR$ defined by
$x_S(\lambda):=\|S(\delta_\lambda)\|_{C_0(E^0)}$
vanishes at infinity on $E^1$.
\end{lemma}

\begin{proof}
First suppose $S=\Theta_{x,y}$ for some $x,y \in X(E)$. Then for $\lambda
\in E^1$, we have
\[
\|\Theta_{x,y}(\delta_\lambda)\|^2
=\sum_{r(\mu)=r(\lambda)} | x(\mu) y(\lambda) |^2
\le |y(\lambda)|^2 \|x\|^2;
\]
since $y \in X(E) \subset C_0(E^1)$, so is $\lambda
\mapsto \|\Theta_{x,y}(\delta_\lambda)\|$.
Easy calculations show that $|x_{wS+zT}(\lambda)|
\leq |w|\,|x_S(\lambda)|+|z|\,|x_T(\lambda)|$ and $|x_S(\lambda)|\leq
\|S\|_{\Ll(X(E))}$, so the result for arbitrary
$S\in
\Kk(X(E))$ follows by linearity and continuity.
\end{proof}

\begin{example}\label{eg:not compactly aligned} (A Cuntz-Krieger bimodule
which is not compactly aligned.) Let $(G,P) = (\ZZ^2,\NN^2)$. Let
$E^0:=\{(0,0), (0,1), (1,0), (1,1)\}$,
\[
E^1_{(1,0)} := \{\lambda\}\cup\{\alpha_i : i \in \NN\},\quad
E^1_{(0,1)} := \{\mu_i : i \in \NN\}\cup\{\beta\},
\]
and define
\begin{gather*}
r(\lambda) = (1,0),\quad s(\lambda) = (0,0),\quad r(\alpha_i)
= (1,1), \quad s(\alpha_i) = (0,1),\quad\text{and}\\
r(\mu_i) = (1,1), \quad s(\mu_i) =
(1,0),\quad r(\beta) = (0,1),\quad s(\beta) = (0,0).
\end{gather*}
By \cite[Theorem~2.1]{FS}, there is a
unique product system $E$ over $\NN^2$ in which $\beta\alpha_i =
\lambda\mu_i$. In
pictures:
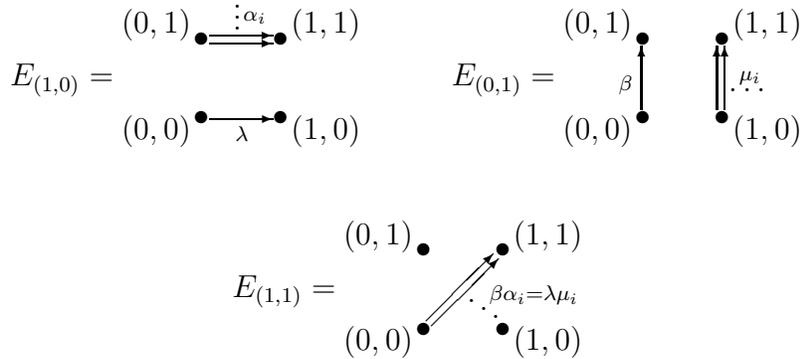
\[
\begin{picture}(300,140)
\putleft(64,105){$E_{(1,0)} = $\phantom{$(0,0)$}}
\putleft(67,85){$(0,0)$}
\putleft(67,125){$(0,1)$}
\put(99,85){\phantom{$\bullet$}$(1,0)$}
\put(99,125){\phantom{$\bullet$}$(1,1)$}
\put(68,90){$\bullet$}
\put(68,120){$\bullet$}
\put(98,90){$\bullet$}
\put(98,120){$\bullet$}
\put(74,92.5){$\vector(1,0){24}$}
\put(74,124){$\vector(1,0){24}$}
\put(74,121){$\vector(1,0){24}$}
\put(84,86){$_{\lambda}$}
\put(83,126.5){$\vdots$}
\put(87,131.5){$_{\alpha_i}$}
\putleft(231,105){$E_{(0,1)} =$\phantom{$(0,0)$}}
\putleft(234,85){$(0,0)$}
\putleft(234,125){$(0,1)$}
\put(266,85){\phantom{$\bullet$}$(1,0)$}
\put(266,125){\phantom{$\bullet$}$(1,1)$}
\put(235,90){$\bullet$}
\put(235,120){$\bullet$}
\put(265,90){$\bullet$}
\put(265,120){$\bullet$}
\put(237.5,96){$\vector(0,1){24}$}
\put(266,96){$\vector(0,1){24}$}
\put(269,96){$\vector(0,1){24}$}
\put(229,105){$_{\beta}$}
\put(270.5,103){$\dots$}
\put(274.5,109){$_{\mu_i}$}

\putleft(148,25){$E_{(1,1)} =$\phantom{$(0,0)$}}
\putleft(151,5){$(0,0)$}
\putleft(151,45){$(0,1)$}
\put(183,5){\phantom{$\bullet$}$(1,0)$}
\put(183,45){\phantom{$\bullet$}$(1,1)$}
\put(152,10){$\bullet$}
\put(152,40){$\bullet$}
\put(182,10){$\bullet$}
\put(182,40){$\bullet$}
\put(156,16){$\vector(1,1){25.5}$}
\put(158,14){$\vector(1,1){25.5}$}
\put(171,16){$\ddots$}
\put(180,25.5){$_{\beta\alpha_i = \lambda\mu_i}$}
\end{picture}
\]
For $S := \Theta_{\delta_{\lambda}, \delta_{\lambda}}$ and $T :=
\Theta_{\delta_{\beta}, \delta_{\beta}}$, we can compute
$S^{(1,1)}_{(1,0)}\circ T^{(1,1)}_{(0,1)}(\delta_{\lambda\mu_i})$ using
(\ref{defSpq}). To evaluate $T^{(1,1)}_{(0,1)}(\delta_{\lambda\mu_i})$ we
need to factor
$\lambda\mu_i$ as $\beta\alpha_i$, so that
$\delta_{\lambda\mu_i}=\delta_\beta\delta_{\alpha_i}$. Then
\begin{align}\label{calculaterank1}
S^{(1,1)}_{(1,0)}\circ T^{(1,1)}_{(0,1)}(\delta_{\lambda\mu_i})
&=S^{(1,1)}_{(1,0)}(T(\delta_\beta)\delta_{\alpha_i})
=S^{(1,1)}_{(1,0)}(\delta_{\beta}\delta_{\alpha_i})\\
&=S^{(1,1)}_{(1,0)}(\delta_\lambda\delta_{\mu_i})
=S(\delta_{\lambda})\delta_{\mu_i}=\delta_{\lambda\mu_i}.\notag
\end{align}
Thus $\lambda\mu_i\mapsto
\big\|S^{(1,1)}_{(1,0)}\circ T^{(1,1)}_{(0,1)}(\delta_{\lambda\mu_i})\big\|$
does not
vanish at infinity  on $E^1_{(1,1)}$. Lemma~\ref{lem:compacts C_0}
therefore implies
that $S^{(1,1)}_{(1,0)}\circ T^{(1,1)}_{(0,1)}$ is not compact, and $E$ is
not compactly aligned.
\end{example}

To identify the $E$ for which $X(E)$ is
compactly aligned, we legislate out the behaviour which makes
Example~\ref{eg:not compactly aligned} work. More precisely:

\begin{definition}\label{def:finitely aligned}
Suppose $(E,\varphi)$ is a
product system of graphs over a quasi-lattice ordered semigroup $P$, and let
$\mu
\in E^1_p$ and
$\nu
\in E^1_q$. A \emph{common extension} of $\mu$ and $\nu$ is a path $\gamma$
such that $\gamma(0,p)=\mu$ and $\gamma(0,q)=\nu$. Notice that  $d(\gamma)$
is then an upper bound for $p$ and $q$, so $p\vee q<\infty$; we say that
$\gamma$ is a {\em minimal common extension} if $d(\gamma)=p\vee q$. We
denote by
$\MCE(\mu,\nu)$ the set of minimal common extensions of $\mu$ and $\nu$, and
say that $(E, \varphi)$ is {\em
finitely aligned} if $\MCE(\mu,\nu)$ is finite (possibly  empty) for
all $\mu,\nu \in
E^1$.
\end{definition}

\begin{theorem}\label{thm:compactly aligned}
Let $(E,\varphi)$ be
a product system of graphs over a quasi-lattice ordered semigroup
$P$. Then $X(E)$ is
compactly aligned if and only if $(E,\varphi)$ is finitely aligned.
\end{theorem}

\begin{proof}
If $\MCE(\lambda,\beta)$ is infinite for some $\alpha$ and $\beta$, there are
infinitely many paths
$\mu_i$ and $\alpha_i$ such that $\lambda\mu_i=\beta\alpha_i$, and the
argument of  Example~\ref{eg:not compactly aligned} shows that $X(E)$ is not
compactly aligned. Suppose that
$(E,\varphi)$ is finitely aligned,
$p,q \in P$ satisfy $p\vee q < \infty$, and $\mu_1, \mu_2 \in E^1_p$,
$\nu_1, \nu_2 \in E^1_q$. Then computations like (\ref{calculaterank1})
show that $(\Theta_{\delta_{\nu_1}, \delta_{\nu_2}})^{p\vee
q}_q(\delta_\lambda)=0$ unless $\lambda(e,q)=\nu_2$,
and then with
$\sigma:=\nu_1\lambda(q,p\vee q)$ we have
\begin{align*}
(\Theta_{\delta_{\mu_1}, \delta_{\mu_2}})^{p\vee q}_p (\Theta_{\delta_{\nu_1},
\delta_{\nu_2}})^{p\vee q}_q(\delta_\lambda)&=
\delta_{\nu_2}(\lambda(0,q))\delta_{\mu_2}(\sigma(0,p))
\delta_{\mu_1\sigma(p,p\vee q)}\\
&=
\begin{cases}\delta_{\mu_1\sigma(p,p\vee q)} &\text{if $\sigma(0,p) =
\mu_2$} \\ 0 &\text{otherwise.} \end{cases}
\end{align*}
Thus
\[
(\Theta_{\delta_{\mu_1}, \delta_{\mu_2}})^{p\vee q}_p (\Theta_{\delta_{\nu_1},
\delta_{\nu_2}})^{p\vee q}_q = \sum_{\sigma \in \MCE(\mu_2, \nu_1)}
\Theta_{\delta_{\mu_1
\sigma(p,p\vee q)}, \delta_{\nu_2 \sigma(q,p\vee q)}},
\]
which belongs to $\Kk(X(E))$ because $\MCE(\mu_2, \nu_1)$ is finite.
\end{proof}

\section{Nica covariance}\label{sec4}
In this section, we show that when $X=X(E)$, Fowler's Nica-covariance condition
reduces to an extra relation for Toeplitz $E$-families, which will
look familiar to
anyone who has studied any generalisation of Cuntz-Krieger algebras.
This relation
automatically holds for Toeplitz-Cuntz-Krieger families of single
graphs,   but is
not automatic for the Toeplitz families of product systems.

Suppose $X$ is a product system of
$A-A$ bimodules over a quasi-lattice ordered semigroup $P$, and $\psi$ is a
nondegenerate Toeplitz representation of $X$ on $\Hh$. Fowler shows in
\cite[Proposition~4.1]{F99} that there is an action
$\alpha^\psi:P\to \End \psi_e(A)^\prime$ such that
\begin{equation}\label{defalpha}
\alpha^\psi_p(T)\psi_p(x)=\psi_p(x)T\mbox{ for }T\in \psi_e(A)^\prime\mbox{
and }\alpha^\psi_p(1)h=0 \mbox{ for }h\in\psi_p(X_p)^\perp.
\end{equation}
The representation
$\psi$ is {\em Nica covariant} if
\begin{equation}\label{eq:alpha N.C.}
\alpha^\psi_p(1_p)\alpha^\psi_q(1_q) =
\begin{cases}
\alpha^\psi_{p\vee q}(1_{p\vee q}) & \text{if $p \vee q < \infty$} \\
0 & \text{otherwise.}
\end{cases}
\end{equation}
We denote by $(\Tt_{\cov}(X),i_X)$ the pair which is universal for
Nica-covariant
Toeplitz representations of $X$ in the sense of
\cite[Theorem~6.3]{F99}. When $X$ is compactly aligned, it follows from
\cite[Lemma~5.5 and Proposition~5.6]{F99} that the Nica covariance condition
\eqref{eq:alpha N.C.} makes sense for a representation taking values in a
$C^*$-algebra, and then $(\Tt_{\cov}(X), i_X)$ is universal in the
usual sense of the
word.

When $P$ is the positive cone in a totally ordered group, $p\vee q$
is either $p$ or
$q$, and Nica covariance is automatic. Thus Toeplitz representations
of a single
Cuntz-Krieger bimodule $X(E)$ are always Nica covariant.  For product
systems of row-finite graphs over lattice-ordered semigroups such as
$\NN^k$, Nica
covariance is a consequence of Cuntz-Pimsner covariance:

\begin{lemma} \label{lem:r.f. -> n.c.}
Let  $(E,\varphi)$ be a
product system of graphs over a lattice-ordered semigroup $P$. If
every $E_p$ is row-finite, then every Toeplitz representation
of
$X(E)$ which is Cuntz-Pimsner covariant is also Nica covariant. In
particular, if
$\Lambda$ is a row-finite $k$-graph, every
Cuntz-Pimsner covariant representation of $X(E_\Lambda)$ is Nica
covariant.
\end{lemma}

\begin{proof}
Since each $E_p$ is row-finite, $C_0(E^0)$ acts by compact
operators on the left of each $X(E_p)$
\cite[Proposition~4.4]{FR}, and the result follows from
\cite[Proposition~5.4]{F99}.
\end{proof}

\begin{cor}
Let  $(E,\varphi)$ be a
product system of row-finite graphs over a lattice-ordered semigroup
$P$. Then $\Oo_{X(E)}$ is isomorphic to a quotient of $\Tt_{\cov}(X(E))$.
\end{cor}

\begin{prop}\label{prop:alpha formula}
Let  $(E,\varphi)$ be a
product system of graphs over a quasi-lattice ordered semigroup $P$,
and let $\psi$ be a nondegenerate Toeplitz representation of
$X(E)$ on $\Hh$. For $p \in P$, $T \in \Bb(\Hh)$ and $h \in \Hh$, the sum
\[
\sum_{\lambda \in E^1_p} \psi_p(\delta_\lambda) T
\psi_p(\delta_\lambda)^* h
\]
converges in $\Hh$; if $T\in \psi_e(C_0(E^0))'$, it converges to
$\alpha_p^\psi(T)h$.
\end{prop}

\begin{proof}
By \cite[Proposition~4.1(1)]{F99}, it suffices to work with a representation
$(\psi,\pi)$ of a single graph
$E$, and show
\begin{itemize}
\item[(1)] that the sum $\alpha(T) h := \sum_{\lambda \in E^1}
\psi(\delta_\lambda)T\psi(\delta_\lambda)^* h$ converges for all $h
\in \Hh$;
\item[(2)] that $\alpha(T) \in \Bb(\Hh)$ for each $T \in \Bb(\Hh)$;
\item[(3)] that $\alpha$ is an endomorphism of $\pi(C_0(E^0))^\prime$; and
\item[(4)] that $\alpha$ satisfies $\alpha(T)\psi(x) = \psi(x)T$ for $T\in
\psi_e(C_0(E^0))'$, and
$\alpha(1)|_{(\psi(X)\Hh)^{\perp}} = 0$.
\end{itemize}

Because the $\psi(\delta_\lambda)$ are
partial isometries with orthogonal ranges, we have
\[
\sum_{\lambda \in E^1} \|\psi(\delta_\lambda)T \psi(\delta_\lambda)^*
h\|^2 \le \sum_{\lambda\in E^1}\|T\|^2\|\psi(\delta_\lambda)^*h\|^2\leq
\|T\|^2\|h\|^2.
\]
Thus $\sum_{\lambda \in E^1}
\psi(\delta_\lambda)T\psi(\delta_\lambda)^*h$ is a sum of orthogonal
vectors which converges
in $\Hh$, and the sum satisfies
\[
\|\alpha(T)h\|^2 = \Big\|\sum_{\lambda\in E^1} \psi(\delta_\lambda)T
\psi(\delta_\lambda)^* h\Big\|^2 =
\sum_{\lambda \in E^1} \|\psi(\delta_\lambda)T \psi(\delta_\lambda)^*
h\|^2\leq \|T\|^2\|h\|^2.
\]
This gives (1) and  (2).

Multiplying $\psi(\delta_\lambda)T\psi(\delta_\lambda)^*$ on either side by
$\psi(\delta_v)$ gives $0$ unless $v=s(\lambda)$, and leaves it alone if
$v=s(\lambda)$. Thus each
$\psi(\delta_\lambda)T\psi(\delta_\lambda)^*$ belongs to
$\pi(C_0(E^0))'$, and so does the strong sum $\alpha(T)$. If $S$ and
$T$ belong to
$\pi(C_0(E^0))'$, then
\begin{align*}
\psi(\delta_\lambda)S\psi(\delta_\lambda)^*\psi(\delta_\mu)T\psi(\delta_\mu)^*
&=\psi(\delta_\lambda)S\psi(\langle\delta_\lambda,\delta_\mu\rangle_{C_0(E^0)})
T\psi(\delta_\mu)^*\\
&=\begin{cases}
\psi(\delta_\lambda)ST\psi(\langle\delta_\lambda,\delta_\mu\rangle_{C_0(E^0)})
\psi(\delta_\mu)^*&\text{if $\mu=\lambda$}\\
0&\text{otherwise}
\end{cases}\\
&=\begin{cases}
\psi(\delta_\lambda)ST\psi(\delta_\lambda)^*&\text{if $\mu=\lambda$}\\
0&\text{otherwise,}
\end{cases}
\end{align*}
and it follows by taking sums and limits that $\alpha$ is multiplicative on
$\pi(C_0(E^0))'$. It is clearly $*$-preserving.

For (4), we let $T\in \psi_e(C_0(E^0))'$ and calculate:
\[
\alpha(T)\psi(\delta_\lambda)=\sum_{\mu\in E^1}
\psi(\delta_\mu)T\psi(\delta_\mu)^*\psi(\delta_\lambda)=
\psi(\delta_\lambda)T\pi(\delta_{r(\lambda)})
=\psi(\delta_\lambda)\pi(\delta_{r(\lambda)})T
=\psi(\delta_\lambda)T.
\]
Extending by linearity gives $\alpha(T)\psi(x) = \psi(x)T$ for $x\in
X_c(E)$, which
suffices by continuity. If
$h\perp\psi(X)\Hh$, then
$\psi(\delta_\lambda)^*h=0$ for all
$\lambda$, and
$\alpha(T)h=0$.
\end{proof}

Suppose that $\{S_\lambda\}\subset B(\Hh)$ is a
Toeplitz $E$-family for a product system $(E,\varphi)$ of graphs
over a quasi-lattice ordered semigroup $P$.
Proposition~\ref{prop:alpha formula} implies that the corresponding Toeplitz
representation $\psi$ of
$X(E)$ is Nica covariant if and only if
\begin{equation} \label{eq:Nica cov}
\textstyle{\Big(\sum_{\mu \in E^1_p}
S_\mu S_\mu^*\Big)
\Big(\sum_{\nu \in E^1_q}
S_\nu S_\nu^*\Big) =
\begin{cases}
{\sum_{\lambda \in E^1_{p\vee q}} S_\lambda S_\lambda^*} &\text{if $p
\vee q < \infty$} \\ 0 &\text{otherwise.}
\end{cases}}
\end{equation}
The sums in \eqref{eq:Nica cov} may be infinite, and then only
converge in the strong operator topology, so this is a spatial
criterion rather than a
$C^*$-algebraic one. When
$E$ is finitely aligned, however, there is an equivalent condition
which only uses
finite sums.

\begin{prop} \label{prop:finite N.C.}
Let  $(E,\varphi)$ be a
finitely aligned product system of graphs over a quasi-lattice ordered
semigroup
$P$, and let $\{S_\lambda\}\subset B(\Hh)$ be a
Toeplitz $E$-family. The corresponding Toeplitz representation
$\psi$ of $X(E)$ is Nica covariant if and only if, for all $p,q \in P$,
$\mu \in E^1_p$ and $\nu \in E^1_q$, we have
\begin{equation} \label{eq:finite N.C.}
S_\mu^*S_\nu
  = \sum_{\mu\alpha=\nu\beta\in\MCE(\mu,\nu)}  S_\alpha
S_\beta^*\quad\text{(which is $0$ if $p\vee q=\infty$).}
\end{equation}
\end{prop}

\begin{proof}
First suppose $\psi$ is Nica covariant, and let $\mu \in E^1_p$ and
$\nu \in E^1_q$. Then because the $S_\lambda$ corresponding to
$\lambda$ of the same
degree have mutually orthogonal ranges, we have
\begin{align*}
S_\mu^*S_\nu&=S_\mu^*\textstyle{\Big(\sum_{\gamma \in
E^1_p} S_\gamma S_\gamma^* \Big)
\Big(\sum_{\sigma \in E^1_q}
S_\sigma S_\sigma^* \Big)}
S_\nu\\
&=
\begin{cases}
   S_\mu^* \big(\sum_{\lambda \in E^1_{p\vee q}} S_\lambda
S_\lambda^*\big) S_\nu
   & \text{if $p\vee q<\infty$} \\
   0
   & \text{if $p\vee q=\infty$}
\end{cases}\\
&=\sum_{\mu\alpha=\nu\beta \;\in\; \MCE(\mu,\nu)} S_\alpha S_\beta^*,
\end{align*}
because $(S_\mu^*S_\lambda)(S_\lambda^*S_\nu)=0$ unless
$\lambda=\mu\alpha=\nu\beta$, and
$\MCE(\mu,\nu)$ is empty if $p\vee q=\infty$.

On the other hand, let $p,q\in P$ and suppose that \eqref{eq:finite
N.C.} holds. Then
\[
\textstyle{\Big(\sum_{\mu \in E^1_p}
S_\mu S_\mu^*\Big)
\Big(\sum_{\nu \in E^1_q}
S_\nu S_\nu^*\Big)}=
\displaystyle{\sum_{\mu \in E^1_p,\nu\in E^1_q}} S_\mu\Big(
\textstyle{\sum_{\mu\alpha=\nu\beta \in\MCE(\mu,\nu)}  S_\alpha
S_\beta^*\Big)S_\nu^*}
\]
which is $\sum\{S_\lambda S_\lambda^*:\lambda\in E^1_{p\vee q}\}$ if
$p \vee q <
\infty$ because the factorisation property implies that each
$\lambda$ appears exactly once as a $\mu\alpha$ and as a $\nu\beta$,
and $0$ if $p
\vee q = \infty$ because then each $\MCE(\mu,\nu)$ is empty.
\end{proof}

\section{Toeplitz-Cuntz-Krieger families}\label{sec5}
Relation (\ref{eq:finite N.C.}) is familiar: some version of it is
used in every
theory of Cuntz-Krieger algebras to ensure that $\lsp\{S_\mu
S_\nu^*\}$ is a dense
$*$-subalgebra of $C^*(\{S_\mu\})$ (see, for example,
\cite[Lemma~2.2]{CK1}, \cite[Lemma~1.1]{KPR}, \cite[Proposition~3.5]{RSY}). As
Lemma~\ref{lem:r.f. -> n.c.} shows, it is often automatic when the graphs are
row-finite, but otherwise it will have to be assumed if we want
$C^*(\{S_\mu\})$
to behave like a Cuntz-Krieger algebra.

We therefore make the following definition:

\begin{definition}\label{def:new TCK family}
Let $E$ be a finitely aligned product system of graphs over a quasi-lattice
ordered semigroup $P$. Partial isometries
$\{s_\lambda :
\lambda
\in E^1\}$ in a
$C^*$-algebra
$B$ form a {\em Toeplitz-Cuntz-Krieger $E$-family} if:
\begin{itemize}
\item[(1)] $\{s_v : v \in E^0\}$ are mutually orthogonal
projections,
\item[(2)] $s_\lambda s_\mu = s_{\lambda\mu}$ for all $\lambda,\mu
\in E^1$ such that $r(\lambda) = s(\mu)$,
\item[(3)] $s^*_\lambda s_\lambda = s_{r(\lambda)}$ for all $\lambda
\in E^1$,
\item[(4)] for all $p \in P\setminus\{e\}, v \in E^0$ and every
finite $F \subset s^{-1}_p(v)$, $s_v \ge \sum_{\lambda
\in F} s_\lambda s^*_\lambda$,
\item[(5)] $s_\mu^*s_\nu
  = \sum_{\mu\alpha=\nu\beta \;\in\; \MCE(\mu,\nu)}  s_\alpha s_\beta^*$ for all
$\mu,\nu \in E^1$.
\end{itemize}
They form a \emph{Cuntz-Pimsner $E$-family} if they also satisfy
\begin{itemize}
\item[(6)] $s_v = \sum_{\lambda \in s_p^{-1}(v)} s_\lambda
s^*_\lambda$ whenever
$s_p^{-1}(v)$ is finite.
\end{itemize}
\end{definition}

\begin{remark}
Multiplying both sides of (5) on the left by $s_\mu$ and on the right
by $s_\nu^*$
gives
\begin{equation}\label{altNicacov}
(s_\mu s_\mu^*)(s_\nu s_\nu^*)=\sum_{\gamma\in \MCE(\mu,\nu)}
s_\gamma s_\gamma^*,
\end{equation}
and this is equivalent to (5) because we can get back by multiplying
on the left by
$s_\mu^*$ and on the right by $s_\nu$.
\end{remark}

\begin{remark}
We have called families satisfying (6) Cuntz-Pimsner families rather than
Cuntz-Krieger families because of the problems with sinks mentioned in
Remark~\ref{bloodysinks}: if $v$ is a sink in a single graph $E$, then (6)
implies that
$s_v=0$, whereas the generally acccepted Cuntz-Krieger relations
impose no relation
at $v$. The Cuntz-Pimsner families are the ones which correspond to
Cuntz-Pimsner covariant representations of $X(E)$.
\end{remark}

\begin{example}[The Fock representation]\label{eg:Fock repn}
For $\lambda\in E^1$, let
$S_\lambda$ be the partial isometry on $\ell^2(E^1)$ such that
\[
S_{\lambda} e_\mu := \begin{cases} e_{\lambda\mu} &\text{if
$r(\lambda) = s(\mu)$} \\ 0 &\text{otherwise.} \end{cases}
\]
We claim that $\{S_\lambda : \lambda \in E^1\}$ is a Toeplitz-Cuntz-Krieger
$E$-family. Conditions (1)--(3) of \defref{def:new TCK family} are
obvious, and (4) holds because
\begin{equation} \label{notCK}
\Big(S_v-\textstyle{\sum_{\lambda \in s_p^{-1}(v)}}S_\lambda
S_\lambda^*\Big)e_v=e_v
\end{equation}
for all $v \in E^0$ and $p \in P\setminus\{e\}$.
To verify (5), we compute on the one hand
\[
\big(S^*_\lambda S_\mu e_\nu | e_\sigma\big)
= \big(S_\mu e_\nu | S_\lambda e_\sigma\big) = \begin{cases}
1 &\text{if $\mu\nu = \lambda\sigma$} \\
0 &\text{otherwise,}
\end{cases}
\]
and on the other hand,
\begin{align*}
\Big( \sum_{\lambda\alpha=\mu\beta \in
\MCE(\lambda,\mu)}S_\alpha
S^*_\beta e_\nu \Big| e_\sigma
\Big)
&= \sum_{\lambda\alpha=\mu\beta \in \MCE(\lambda,\mu)}
\big(S^*_\beta e_\nu |
S^*_\alpha e_\sigma\big) \\
&= \sum_{\lambda\alpha=\mu\beta \in \MCE(\lambda,\mu)}
\begin{cases}
1 &\text{if $\nu = \beta\tau$ and $\sigma = \alpha\tau$ for some $\tau$} \\
0 &\text{otherwise.}
\end{cases}
\end{align*}
By the factorisation property, at most one term in this
last sum can be nonzero, and there is one precisely when
$\lambda\alpha\tau = \mu\beta\tau$ for some $\lambda\alpha=\mu\beta \in
\MCE(\lambda,\mu)$, giving (5).

If there is a vertex $v$ which emits just finitely many edges in
some $E_p$, then \eqref{notCK} implies that (6) does not hold, and hence
$\{S_\lambda\}$ is not a Cuntz-Pimsner family.
\end{example}

If $(E, \varphi)$ is finitely aligned, then Theorem~\ref{thm:families
and repns} and
\propref{prop:finite N.C.} imply that the Toeplitz $E$-family
$\{i_{X(E)}(\delta_\lambda) : \lambda \in E^1\}$ in $\Tt_{\cov}(X(E))$ is a
Toeplitz-Cuntz-Krieger $E$-family. It then follows from
\lemref{lem:xcdense} that
$\Tt_{\cov}(X(E))$ is generated by $\{i_{X(E)}(\delta_\lambda)\}$. We
can now apply
the other direction of \thmref{thm:families and repns} to see that
$\Tt_{\cov}(X(E))$
is universal for Toeplitz-Cuntz-Krieger $E$-families. Thus:

\begin{cor}\label{cor:Tcov univ}
Let $(E,\varphi)$ be a finitely aligned product system of graphs over
a quasi-lattice
ordered semigroup $P$. Then $(\Tt_{\cov}(X(E)),
\{i_{X(E)}(\delta_\lambda)\})$ is
universal for Toeplitz-Cuntz-Krieger $E$-families.
\end{cor}

In view of \corref{cor:Tcov univ}, we define $\Tt C^*(E)$ to be the
universal algebra
$\Tt_{\cov}(X(E))$. If there are no sinks, we define $C^*(E)$ to be
the quotient of
$\Tt C^*(E)$ which is universal for Cuntz-Pimsner $E$-families. If
$\Lambda$ is a
row-finite
$k$-graph with no sources, it follows from \lemref{lem:r.f. -> n.c.} that
$C^*(E_\Lambda)$ is the
$C^*$-algebra
$C^*(\Lambda)$ studied in \cite{KP}.

  From now on, we denote by
$\{s_\lambda:\lambda\in E^1\}$ the canonical generating family in
$\Tt C^*(E)$, and
if $\{t_\lambda:\lambda\in E^1\}$ is a Toeplitz-Cuntz-Krieger $E$-family in a
$C^*$-algebra $B$, then we write
$\pi_t$ for the homomorphism of $\Tt C^*(E)$ into $B$ such that
$\pi_t(s_\lambda)=t_\lambda$.

We now see what Fowler's theory tells us about faithful representations.

\begin{prop} \label{thm:7.2}
Let $(G,P)$ be quasi-lattice ordered with $G$
amenable, and let $(E,\varphi)$ be a finitely aligned product system of
graphs over $P$. Let $\{S_\lambda : \lambda \in E^1\}$ be
a Toeplitz-Cuntz-Krieger $E$-family in $\Bb(\Hh)$, and suppose that
for every finite subset $R$ of $P\setminus\{e\}$ and every $v \in E^0$, we have
\begin{equation} \label{eq:7.2}
\prod_{p\in R}\Big(S_v  - \textstyle{\sum_{\lambda \in
s^{-1}_{p}(v)}} S_\lambda
S^*_\lambda\Big) > 0.
\end{equation}
Then the corresponding representation $\pi_S:\Tt C^*(E)\to B(\Hh)$ is faithful.
\end{prop}

\begin{proof}
We consider the representation $\psi$ of $X(E)$ associated to $\{S_\lambda\}$.
Theorem~\ref{thm:compactly aligned} says that $X(E)$ is compactly aligned, and
Proposition~\ref{prop:finite N.C.} that $\psi$ is Nica covariant. Since the
$\delta_v$ span a dense subspace of
$C_0(E^0)$ and the $\psi_e(\delta_v)=S_v$ are mutually orthogonal,
\propref{prop:alpha formula} implies that \eqref{eq:7.2} is
equivalent to the displayed hypothesis in \cite[Theorem~7.2]{F99}. Thus
\cite[Theorem~7.2]{F99} implies that $\psi_*$ is faithful on
$\Tt_{\cov}(X(E))$. But
$\pi_S$ is by definition the representation $\psi_*$ of $\Tt
C^*(E):=\Tt_{\cov}(X(E))$.
\end{proof}

\begin{cor} \label{cor:Fock faithful}
Let $(G,P)$ be a quasi-lattice ordered group such that $G$ is
amenable, and let $(E,\varphi)$ be a finitely aligned product system of
graphs over $P$. Then the representation $\pi_S$ of $\Tt C^*(E)$
associated to the
Fock representation of Example~\ref{eg:Fock repn} is faithful.
\end{cor}
\begin{proof}
Equation~(\ref{eq:7.2}) follows from (\ref{notCK}).
\end{proof}

\section{A $C^*$-algebraic uniqueness theorem}\label{sec6}

\begin{theorem} \label{thm:faithfulness}
Let $(G,P)$ be a quasi-lattice ordered
group such that $G$ is amenable, and  let $(E,\varphi)$ be a finitely
aligned product system of graphs over $P$. Let $\{t_\lambda
: \lambda \in E^1\}$ be a Toeplitz-Cuntz-Krieger $E$-family in a
$C^*$-algebra $B$.
Suppose that for every finite subset $R$ of
$P\setminus\{e\}$, every $v \in E^0$, and every
collection of finite sets $F_p \subset s^{-1}_{p}(v)$, we have
\begin{equation}\label{eq:faithfulness}
\prod_{p\in R} \Big(t_v - \textstyle{\sum_{\lambda \in F_p} t_\lambda
t^*_\lambda}\Big) > 0.
\end{equation}
Then the associated homomorphism $\pi_t:\Tt C^*(E)\to B$ is injective.
\end{theorem}

To prove \thmref{thm:faithfulness}, we first establish that there is
a linear map
$\Phi^E$ onto the diagonal in  $\Tt C^*(E)$ which is faithful on positive
elements, and show that there is a norm-decreasing linear map
$\Phi^B$ on $\pi_t(\Tt C^*(E))$ such that $\pi_t \circ \Phi^E =
\Phi^B \circ \pi_t$.

\begin{prop} \label{prop:faithful expectation}
There is a linear map $\Phi^E : \Tt C^*(E) \to
\Tt C^*(E)$ such that
\[
\Phi^E(s_\lambda s_\mu^*) =
\begin{cases}
s_\lambda s_\lambda^* &\text{if $\lambda = \mu$}
\\
0 &\text{otherwise,}
\end{cases}
\]
and $\Phi^E$ is faithful on positive elements.
\end{prop}

\begin{proof}
Let $\{e_i : i \in I\}$ be an orthonormal basis for $\Hh$, and for $i
\in I$, let $P_i$ be the projection onto $\CC e_i$. Then for $T\in B(\Hh)$,
$\sum_{i\in I} P_i T P_i$ converges in the strong operator
topology, and $T \mapsto \sum_{i\in I} P_i T P_i$ is the diagonal map
on $\Bb(\Hh)$ which takes  the rank-one operator $\Theta_{e_i, e_j}$
to $\Theta_{e_i,
e_i}$ if $i = j$ and to $0$ otherwise. It follows that this diagonal
map is linear and norm-decreasing, and  it is
faithful on positive elements: $\Phi(T^*T)=0$ implies
$(T^*Te_i|e_i)=0$ for all $i$,
and hence $T=0$.

Let $\Hh := \ell^2(E^1)$ and let $\{S_\lambda : \lambda \in E^1\}$ be
the Toeplitz-Cuntz-Krieger family of Example~\ref{eg:Fock repn}. Then
a calculation using the basis elements $\{e_\nu : \nu \in
E^1\}$ shows that
\[
P_\gamma S_\lambda S^*_\mu P_\gamma =
\begin{cases}
P_\gamma &\text{if $\lambda = \mu = \gamma(e, d(\mu))$} \\
0 &\text{otherwise.}
\end{cases}
\]
Thus if $\Phi$ denotes the diagonal map on $\ell^2(E^1)$, then
\[
\Phi(S_\lambda S^*_\mu) = P_{\clsp\{e_\gamma : \lambda = \mu =
\gamma(e, d(\mu))\}} =
\begin{cases}
S_\lambda S^*_\lambda &\text{if $\lambda = \mu$} \\
0 &\text{otherwise.}
\end{cases}
\]
Because the  representation $\pi_S$ associated to the Fock
representation is faithful
by \corref{cor:Fock faithful}, and because $\Phi$ has the required
properties, we can
pull $\Phi$ back to $\Tt C^*(E)$ to get the required map $\Phi^E$.
\end{proof}

We must now establish the existence of $\Phi^B : \pi_t(\Tt C^*(E)) \to
\pi_t(\Tt C^*(E))$ and show that $\pi_t$ is faithful on
$\Phi^E(\Tt C^*(E))$. To do this, we analyse the structure
of the diagonal $\Phi^E(\Tt C^*(E))$. Since $\Tt C^*(E)$ is spanned
by elements of the form $s_\lambda s_\mu^*$, we
consider the image of $\lsp\{s_\lambda s_\mu^*: \lambda, \mu \in E^1\}$ in the
diagonal. We show that for a finite subset $F$ of $E^1$, $C^*(\{t_\lambda
t^*_\lambda : \lambda \in F\})$ sits inside a finite-dimensional
diagonal subalgebra of $B$, and use the matrix units in
this diagonal subalgebra to show that $\Phi^B$ exists and is
norm-decreasing. We can
then show that $\pi_t$ is faithful on $\lsp\{s_\lambda s_\lambda^*: \lambda\in
E^1\}$ just by checking that the
matrix units are nonzero.

Condition (5) of Definition~\ref{def:new TCK family} shows that
$C^*(\{t_\lambda
t^*_\lambda : \lambda \in F\})$ is typically bigger than
$\lsp\{t_\lambda
t^*_\lambda : \lambda \in F\}$; the two can only be equal if
$\lambda,\mu \in F$ implies $\MCE(\lambda,\mu) \subset F$. Thus we
need to pass to a larger finite set $H$ such that $\lambda,\mu\in H$ imply
$\MCE(\lambda,\mu) \subset H$.

\begin{definition} \label{def:vee F}
For each finite subset $F$ of $E^1$, let
\[
\MCE(F) := \{\lambda \in E^1 : d(\lambda) = \textstyle{\bigvee_{\alpha
\in F}}d(\alpha) \text{ and } \lambda(e, d(\alpha)) = \alpha\text{ for all }
\alpha \in F\},
\]
and let $\vee F := \bigcup_{G \subset F} \MCE(G)$.
\end{definition}

\defref{def:vee F} is consistent with \defref{def:finitely aligned}, since
$\MCE(\{\lambda,\mu\}) =
\MCE(\lambda,\mu)$.

\begin{lemma} \label{lem:vee F}
Let $F$ be a finite subset of $E^1$. Then
\begin{itemize}
\item[(1)] $F \subset \vee F$;
\item[(2)] $\vee F$ is the union of the disjoint sets $\vee
\{\lambda
\in F : s(\lambda) = v\}$ over $v\in s(F)$;
\item[(3)] $\vee F$ is finite; and
\item[(4)] $G \subset \vee F$ implies $\MCE(G) \subset \vee F$.
\end{itemize}
\end{lemma}
\begin{proof}
(1) For $\lambda \in F$, $\{\lambda\} \subset F$ and $\lambda \in
\MCE(\{\lambda\})$.

(2) If $\lambda,\mu\in G$ and $s(\lambda)\not= s(\mu)$, then $\MCE
(G)$ is empty.

(3) It suffices to show that if $F \subset E^1$ is finite, then
$\MCE(F)$ is finite. When $|F| = 1$, this assertion is trivial.
Suppose as an inductive hypothesis that $\MCE(F)$ is finite whenever
$|F| \le k$ for some $k \ge 1$, and suppose that $|F| = k+1$. Let
$\lambda \in F$,
and let $F' := F \setminus\{\lambda\}$. Suppose that $\gamma \in
\MCE(F)$. Since $\gamma(e, \bigvee_{\alpha \in F'} d(\alpha)) \in
\MCE(F')$, we have $\gamma \in \MCE(\lambda,\mu)$ for some $\mu \in
\MCE(F')$. Hence $|\MCE(F)| \le \sum_{\mu \in \MCE(F')} |\MCE(\lambda,\mu)|$.
Each term in this sum is finite because $(E,\varphi)$ is finitely aligned,
and the sum has only finitely many terms by the inductive hypothesis. Hence
$\MCE(F)$ is finite.

(4) Let $G \subset \vee F$ and for $\alpha \in G$ choose $G_\alpha \subset F$
such that $\alpha \in \MCE(G_\alpha)$. Let $H := \bigcup_{\alpha \in G}
G_\alpha$. We will show that $\MCE(G) \subset \MCE(H) \subset
\vee F$. Suppose $\lambda \in \MCE(G)$. Then $d(\lambda) =
\bigvee_{\alpha \in G} d(\alpha) = \bigvee_{\alpha \in G}
\big(\bigvee_{\beta \in G_\alpha} d(\beta)\big) = \bigvee_{\beta \in
H} d(\beta)$. For $\beta \in H$, choose $\alpha \in G$ such that
$\beta \in G_\alpha$. Then $\lambda(e, d(\beta)) = \alpha(e, d(\beta))
= \beta$. Thus $\lambda \in \MCE(H)$.
\end{proof}

It follows from \lemref{lem:vee F}(4) that $\lambda, \mu
\in \vee F$ implies that $\MCE(\lambda,\mu) \subset \vee F$.
Consequently, \lemref{lem:vee F}(1) and (\ref{altNicacov}) imply that
\[
C^*(\{t_\lambda t^*_\lambda : \lambda \in F\}) \subset
C^*(\{t_\lambda t^*_\lambda : \lambda \in \vee F\}) =
\lsp\{t_\lambda t^*_\lambda : \lambda \in \vee F\}.
\]
To write this as a diagonal matrix algebra, we need to be able to
orthogonalise the range projections associated to the edges in $\vee
F$.

\begin{lemma} \label{lem:Qs nonzero}
Let $\lambda \in E^1$. If $F\subset s^{-1}(r(\lambda))$ is finite and
$r(\lambda) \not\in F$, then
\[
t_\lambda t^*_\lambda \Big(\prod_{\mu \in F} (t_{s(\lambda)} -
t_{\lambda\mu} t^*_{\lambda\mu})\Big) > 0.
\]
\end{lemma}
\begin{proof}
We have
\[
\Big\|t_\lambda t^*_\lambda \Big(\prod_{\mu \in F} (t_{s(\lambda)} -
t_{\lambda\mu} t^*_{\lambda\mu})\Big)\Big\|
= \Big\|\prod_{\mu \in F} (t_\lambda t^*_\lambda - t_{\lambda\mu}
t^*_{\lambda\mu})\Big\|
= \Big\|t_\lambda \Big(\prod_{\mu \in F} (t_{r(\lambda)} - t_\mu
t^*_\mu) \Big) t^*_\lambda\Big\|,
\]
which is nonzero by \eqref{eq:faithfulness}.
\end{proof}

We now define our matrix units. First note that (\ref{altNicacov}) for
the Toeplitz-Cuntz-Krieger family $\{t_\lambda\}$ implies that the
range projections
$t_\lambda t_\lambda^*$ commute with each other. Thus for every
finite subset $F$ of
$E^1$ and every
$\lambda \in \vee F$, the operator $Q^{\vee F}_{\lambda}$ defined by
\[
Q^{\vee F}_{\lambda} := t_\lambda t^*_\lambda
\Big(\prod_{\lambda\alpha \in \vee F,\, d(\alpha) \not=
e}(t_{s(\lambda)} - t_{\lambda\alpha} t^*_{\lambda\alpha})\Big)
\]
is a projection which commutes with every $t_\mu t_\mu^*$.

\begin{prop} \label{prop:partition}
Let $F$ be a finite subset of $E^1$ such that $\lambda \in F$ implies
$s(\lambda) \in F$. Then $\{Q^{\vee F}_\lambda : \lambda \in \vee F\}$
is a collection of nonzero mutually orthogonal projections in $B$ such
that $\lsp\{Q^{\vee F}_\lambda : \lambda \in \vee F\} =
\lsp\{t_\lambda t^*_\lambda : \lambda \in \vee F\}$. In particular,
\begin{equation}\label{eq:partition}
\sum_{\lambda \in \vee F} Q^{\vee F}_\lambda = \sum_{v \in s(F)} t_v.
\end{equation}
\end{prop}

The key to proving \propref{prop:partition} is establishing
\eqref{eq:partition}, which we do by induction on $\lvert F
\rvert$. This requires two technical lemmas.

\begin{lemma} \label{lem:maximal subpaths}
Let $F$ be as in \propref{prop:partition}, suppose $\lambda \in
F\setminus E^0$ and let $G := F \setminus
\{\lambda\}$. Then for every $\gamma \in \vee F \setminus \vee G$
there is a unique $\mu_\gamma \in \vee G$ such that
\begin{equation} \label{eq:maximal subpath}
\text{if $\mu \in \vee G$ and  $\gamma(e, d(\mu)) = \mu$ then $d(\mu)
\le d(\mu_\gamma)$.}
\end{equation}
We then have $\gamma \in \MCE(\mu_\gamma, \lambda)$; in particular,
$d(\gamma) = d(\mu_\gamma) \vee d(\lambda)$.
\end{lemma}

\begin{proof}
For $\gamma \in \vee F \setminus \vee G$, let $(\vee G)_\gamma :=
\{\mu \in \vee G : \gamma(e, d(\mu)) = \mu\}$, which is
nonempty because $s(\gamma) \in (\vee G)_\gamma$. For every $\mu \in
(\vee G)_\gamma$, $d(\mu) \le d(\gamma)$, so $d := \bigvee_{\mu \in
(\vee G)_\gamma} d(\mu)$ satisfies $d \le d(\gamma)$.
\lemref{lem:vee F}(4) shows that $\gamma(e, d) \in \vee G$, and then
$\mu_\gamma :=
\gamma(e,d)$ has the required property. To see that $\gamma \in
\MCE(\mu_\gamma, \lambda)$, notice that $\gamma \in \vee F \setminus
\vee G$ implies $\gamma \in \MCE(\mu,\lambda)$ for some $\mu \in \vee
G$. Thus $\mu \in (\vee G)_\gamma$, $d(\mu) \le d(\mu_\gamma)$, and
\[
d(\gamma) = d(\mu) \vee d(\lambda) \le
d(\mu_\gamma) \vee d(\lambda).
\]
On the other hand, we have $d(\gamma) \ge
d(\mu_\gamma)$ by definition, and $d(\gamma) \ge d(\lambda)$ since
$\gamma \in \MCE(\lambda,\mu)$. Hence $d(\gamma) = d(\mu_\gamma)
\vee d(\lambda)$, and $\gamma \in \MCE(\mu_\gamma, \lambda)$.
\end{proof}

\begin{lemma}\label{lem:projections equal}
Let $F$ be as in \propref{prop:partition}, suppose $\lambda \in F
\setminus E^0$ and let $G := F \setminus \{\lambda\}$. Then for
each $\delta \in \vee F \setminus \vee G$,
\begin{equation} \label{eq:QF and QG}
Q^{\vee F}_\delta = Q^{\vee G}_{\mu_\delta} t_\delta t^*_\delta.
\end{equation}
\end{lemma}
\begin{proof}
We shall show that
\begin{itemize}
\item[(1)] $Q^{\vee F}_\delta = Q^{\vee G}_{\mu_\delta} Q^{\vee
F}_\delta$, and
\item[(2)] $Q^{\vee G}_{\mu_\delta} t_{\delta\varepsilon}
t^*_{\delta\varepsilon} = 0$ whenever $\delta\varepsilon \in \vee F$
and $d(\varepsilon) \not=e$,
\end{itemize}
and then use these to prove \eqref{eq:QF and QG}.

To prove (1), let $\delta \in \vee F \setminus \vee G$. Since
$t_{\mu_\delta} t^*_{\mu_\delta} \geq t_\delta t^*_\delta$,
\[
Q^{\vee G}_{\mu_\delta} Q^{\vee F}_\delta
= t_\delta t^*_\delta \Big(\prod_{\mu_\delta\nu \in \vee G,\, d(\nu) \not=
e}
(t_{s(\delta)} - t_{\mu_\delta\nu} t^*_{\mu_\delta\nu})\Big) Q^{\vee
F}_\delta.
\]
Suppose $\mu_\delta\nu \in \vee G$ and $d(\nu) \not= e$. Then
\[
t_\delta t^*_\delta (t_{s(\delta)} - t_{\mu_\delta\nu}
t^*_{\mu_\delta\nu})
= t_\delta t^*_\delta - \sum_{\gamma \in \MCE(\delta,\mu_\delta\nu)}
t_\gamma t^*_\gamma \quad\text{by \eqref{altNicacov}}.
\]
Now suppose $\gamma \in \MCE(\delta, \mu_\delta\nu)$. Then $d(\mu_\gamma) \ge
d(\mu_\delta \nu)$ because $\mu_\delta\nu \in \vee G$, and
$d(\mu_\delta \nu) > d(\mu_\delta)$ because $d(\nu) \not= e$.
In particular $\gamma \not= \delta$. But $\gamma(e,
d(\delta)) = \delta$ because $\gamma \in \MCE(\delta, \mu_\delta
\nu)$. Hence there exists $\varepsilon \in E^1$ such that
$d(\varepsilon) \not= e$ and $\gamma = \delta\varepsilon$. Since $\delta$
and $\mu_\delta \nu$ are in $\vee F$, \lemref{lem:vee
F}(4) ensures that $\gamma \in \vee F$, so $t_{s(\delta)} -
t_\gamma t^*_\gamma$ is a factor in $Q^{\vee F}_\delta$, and
$t_\gamma t^*_\gamma Q^{\vee F}_\delta = 0$. Thus
\[
t_\delta t^*_\delta (t_{s(\delta)} - t_{\mu_\delta\nu}
t^*_{\mu_\delta\nu}) Q^{\vee F}_\delta = t_\delta t^*_\delta Q^{\vee
F}_\delta -\Big(\sum_{\gamma \in \MCE(\delta,\mu_\delta\nu)}
t_\gamma t^*_\gamma\Big)Q^{\vee F}_\delta= Q^{\vee F}_\delta.
\]
Applying this equation to each $\mu_\delta \nu \in \vee G$ with
$d(\nu) \not= e$
establishes (1).

To prove (2), suppose that $\delta\varepsilon \in \vee F$ with
$d(\varepsilon) \not= e$. Then $\mu_{\delta\varepsilon} \in \vee G$, and
$\mu_{\delta\varepsilon}\not=\mu_\delta$: if $\mu_{\delta\varepsilon}
= \mu_\delta$,
then
$d(\delta\varepsilon) = d(\lambda) \vee d(\mu_{\delta\varepsilon}) =
d(\lambda) \vee
d(\mu_\delta) = d(\delta)$, contradicting $d(\varepsilon) \not= e$. However,
$(\delta\varepsilon)(e, d(\mu_\delta)) = \delta(e, d(\mu_\delta)) =
\mu_\delta$, so \lemref{lem:maximal subpaths} implies that
$d(\mu_\delta) <
d(\mu_{\delta\varepsilon})$, and $\mu_{\delta\varepsilon} =
\mu_\delta\alpha$ for some $\alpha$ with $d(\alpha) \not= e$. Since
$\mu_{\delta\varepsilon} \in \vee G$, it follows that
\[
Q^{\vee G}_{\mu_\delta} t_{\delta\varepsilon} t^*_{\delta\varepsilon}
\le (t_{s(\mu_\delta)} - t_{\mu_\delta \alpha} t^*_{\mu_\delta
\alpha}) t_{\delta\varepsilon} t^*_{\delta\varepsilon},
\]
which vanishes because $\mu_\delta \alpha =
(\delta\varepsilon)(e, d(\mu_{\delta\varepsilon}))$. This gives (2).

To finish off, we compute:
\begin{flalign*}
&&Q^{\vee F}_\delta
&= Q^{\vee G}_{\mu_\delta} Q^{\vee F}_\delta \quad\text{by (1)}&\\
&&&= Q^{\vee G}_{\mu_\delta} \Big(\prod_{\delta\varepsilon
\in \vee F,\, d(\varepsilon) \not= e} (t_{s(\mu_\delta)} -
t_{\delta\varepsilon} t^*_{\delta\varepsilon}) \Big) t_\delta
t^*_\delta &\\
&&&= Q^{\vee G}_{\mu_\delta} t_\delta t^*_\delta \quad\text{by
(2).}&\qed
\end{flalign*}
\renewcommand{\qed}{}\end{proof}

\begin{proof}[Proof of \propref{prop:partition}]
The $Q^{\vee F}_\lambda$ are nonzero by \lemref{lem:Qs nonzero}. To
see that the $Q^{\vee F}_\lambda$ are orthogonal, suppose that
$\lambda \not = \mu \in \vee F$. If $d(\lambda) = d(\mu)$ then
$Q^{\vee F}_\lambda Q^{\vee F}_\mu \le t_\lambda t^*_\lambda t_\mu
t^*_\mu = 0$ by (4) of \defref{def:new TCK family}. So suppose that
$d(\lambda) \not=
d(\mu)$. We can assume without loss of generality that
$d(\lambda) \vee d(\mu) > d(\lambda)$. Then $\gamma \in
\MCE(\lambda,\mu)$ implies $\gamma = \lambda\alpha$ where $d(\alpha) \not=
e$, and $\gamma \in \vee F$ by \lemref{lem:vee F}(4). Thus
\eqref{altNicacov} shows that
\[
Q^{\vee F}_\lambda Q^{\vee F}_\mu
\le \Big(\sum_{\gamma \in \MCE(\lambda,\mu)} t_\gamma t^*_\gamma\Big) Q^{\vee
F}_\lambda = 0.
\]

Assuming that \eqref{eq:partition} has been established, let $\lambda
\in \vee F$ and calculate:
\begin{align}
t_\lambda t^*_\lambda
&= t_\lambda t^*_\lambda \Big(\sum_{\mu \in \vee F} Q^{\vee
F}_\mu\Big) \quad\text{by \eqref{eq:partition}} \notag\\
&= \sum_{\mu \in \vee F} \Big(t_\lambda t^*_\lambda t_\mu t^*_\mu
\Big(\prod_{\mu\alpha \in \vee F,\, d(\alpha) \not=
e}(t_{s(\mu)} - t_{\mu\alpha} t^*_{\mu\alpha})\Big)\Big) \notag\\
&= \sum_{\mu \in \vee F} \Big(\Big(\sum_{\gamma \in \MCE(\lambda,
\mu)} t_\gamma t^*_\gamma\Big) \Big(\prod_{\mu\alpha
\in \vee F,\, d(\alpha) \not= e}(t_{s(\mu)} - t_{\mu\alpha}
t^*_{\mu\alpha})\Big)\Big).\label{dag}
\end{align}
Suppose $\mu \in \vee F$ and $\mu \not= \lambda\lambda'$ for any path
$\lambda'$, and that $\gamma \in \MCE(\lambda, \mu)$. \lemref{lem:vee
F}(4) ensures that $\gamma \in \vee F$, and $\gamma \not=\mu$
because $\mu \not= \lambda\lambda'$. Thus $\gamma = \mu\alpha$ for
some path $\alpha$ such that $d(\alpha) \not= e$. Hence the product in
$(\ref{dag})$ vanishes for such
$\mu$, and $(\ref{dag})$ collapses to
\[
t_\lambda t^*_\lambda = \sum_{\lambda\lambda^\prime \in \vee F}
Q^{\vee F}_{\lambda\lambda^\prime}.
\]

It therefore suffices to establish \eqref{eq:partition}. Indeed, $Q^{\vee
F}_\lambda \le s(\lambda)$ for all $\lambda$, so \lemref{lem:vee
F}(2) shows that it suffices to establish \eqref{eq:partition} when $F
\subset s^{-1}(v)$ for some $v \in E^0$. We do this by induction on
$\lvert F \rvert$. Recall that $\lambda \in F$ implies $s(\lambda)
\in F$, so if $\lvert F \rvert = 1$ then $F =\vee F =
\{v\}$ and $Q^{\vee F}_v = t_v$.

Suppose that $\lvert F \rvert = k+1 \ge 2$, and that the proposition
holds for all subsets of $s^{-1}(v)$ containing $v$ and having at most
$k$ elements. Since $\lvert F \rvert > 1$ there exists $\lambda \not=
v$ in $F$. Let $G := F \setminus \{\lambda\}$. For $\mu \in \vee
G$, we have
\[
Q^{\vee F}_\mu
= t_\mu t^*_\mu
\Big(\prod_{\mu\alpha \in \vee G,\, d(\alpha) \not= e}
\big(t_v - t_{\mu\alpha} t^*_{\mu\alpha}\big)\Big)
\Big(\prod_{\gamma = \mu\beta \in \vee F \setminus \vee G}
\big(t_v - t_\gamma t^*_\gamma \big)\Big).
\]
Suppose that $t_v - t_\gamma t^*_\gamma$ is a factor in the second
product and  $\mu_\gamma \not= \mu$. Then $\mu_\gamma = \mu\alpha$ for
some $\alpha$  such that $d(\alpha) \not= e$ because $\mu_\gamma$ is the
maximal subpath of $\gamma$ in $\vee G$. Thus $t_v - t_\gamma
t^*_\gamma$ is larger than the factor $t_v - t_{\mu_\gamma}
t^*_{\mu_\gamma}$ from the first product. So such terms in the second
product can be deleted without changing the product, and we have
\[
Q^{\vee F}_\mu = Q^{\vee G}_\mu \Big(\prod_{\gamma \in \vee F
\setminus \vee G,\, \mu_\gamma = \mu}(t_v - t_\gamma
t^*_\gamma)\Big).
\]
Thus
\begin{align*}
\sum_{\lambda \in \vee F} Q^{\vee F}_\lambda
&= \sum_{\mu \in \vee G} Q^{\vee G}_\mu \Big(\prod_{\gamma
\in \vee F \setminus \vee G,\, \mu_\gamma = \mu}(t_v - t_\gamma
t^*_\gamma)\Big) + \sum_{\delta \in \vee F \setminus \vee G} Q^{\vee
F}_\delta \,\\
&= \sum_{\mu \in \vee G} \Big(Q^{\vee G}_\mu
\Big(\prod_{\gamma \in \vee F \setminus \vee G,\,
\mu_\gamma = \mu}(t_v - t_\gamma t^*_\gamma)\Big) +
\sum_{\delta \in \vee F \setminus \vee G,\, \mu_\delta =
\mu} Q^{\vee F}_\delta \Big)
\end{align*}
by \lemref{lem:maximal subpaths}, and \lemref{lem:projections equal} gives
\begin{align}
\sum_{\lambda \in \vee F} Q^{\vee F}_\lambda
&= \sum_{\mu \in \vee G}\Big(Q^{\vee G}_\mu
\Big(\prod_{\gamma \in \vee F \setminus \vee G,\,
\mu_\gamma = \mu}(t_v - t_\gamma t^*_\gamma)\Big) +
\sum_{\delta \in \vee F \setminus \vee G,\, \mu_\delta =
\mu} Q^{\vee G}_{\mu_\delta} t_\delta t^*_\delta \Big)\notag \\
&= \sum_{\mu \in \vee G} Q^{\vee G}_\mu
\Big(\Big(\prod_{\gamma \in \vee F \setminus \vee G,\,
\mu_\gamma = \mu}(t_v - t_\gamma t^*_\gamma)\Big) +
\sum_{\delta \in \vee F \setminus \vee G,\, \mu_\delta =
\mu} t_\delta t^*_\delta \Big).\label{lastone}
\end{align}
If $\mu \in \vee G$ and $\delta \in \vee F \setminus \vee G$
satisfies $\mu_\delta = \mu$, then \lemref{lem:maximal subpaths}
implies that $d(\delta) = d(\mu) \vee d(\lambda)$. Thus $\{t_\delta
t^*_\delta:\mu_\delta = \mu\}$ are mutually orthogonal, and
(\ref{lastone}) is just
$\sum_{\mu \in \vee G}
Q^{\vee G}_\mu$.
Applying the inductive hypothesis to $G$ now establishes
\eqref{eq:partition} for the given $F$.
\end{proof}

\begin{prop}\label{prop:norm decreasing}
There is a norm-decreasing linear map
\[
\Phi^B : C^*(\{t_\lambda : \lambda \in E^1\}) \to \clsp\{t_\lambda
t^*_\lambda : \lambda \in E^1\}
\]
such that $\Phi^B \circ \pi_t = \pi_t \circ \Phi^E$.
\end{prop}
\begin{proof}
It suffices to show that if $F \subset E^1$ is
finite and $\{\alpha_{\lambda,\mu} : \lambda, \mu \in F\} \subset
\CC$, then $\big\|\sum_{\lambda,\mu \in F}
\alpha_{\lambda,\mu} t_\lambda t^*_\mu
\big\| \ge \big\|\sum_{\lambda \in F} \alpha_{\lambda,\lambda}
t_\lambda t^*_\lambda\big\|$.

Since $\sum_{\gamma \in F} Q^{\vee F}_{\gamma}
= \sum_{v \in s(F)} t_v$ and the $Q^{\vee F}_{\gamma}$ commute with
the $t_\lambda
t_\lambda^*$, there exists $\gamma\in\vee F$ such that
\begin{equation}\label{eq:norm attained}
\Big\| Q^{\vee F}_\gamma \Big(\sum_{\lambda \in F}
\alpha_{\lambda,\lambda} t_\lambda t^*_\lambda\Big) \Big\|
= \Big\|\sum_{\lambda \in F} \alpha_{\lambda,\lambda} t_\lambda
t^*_\lambda \Big\|.
\end{equation}
If $\lambda \in F$ and $\gamma \not= \lambda\beta$ for any
$\beta$, then $\delta \in \MCE(\lambda,\gamma)$ implies
$d(\delta) > d(\gamma)$, giving
\[
Q^{\vee F}_\gamma t_\lambda = Q^{\vee F}_\gamma t_\lambda t^*_\lambda
t_\lambda =
\Big(\prod_{\gamma\beta \in \vee F,\, d(\beta) \not= e}
\big(t_\gamma t^*_\gamma - t_{\gamma\beta} t^*_{\gamma\beta}\big)\Big)\Big(
\sum_{\delta \in \MCE(\gamma, \lambda)}
t_\delta t^*_\delta\Big) t_\lambda = 0.
\]
Thus
\[
Q^{\vee F}_\gamma \Big(\sum_{\lambda,\mu \in F} \alpha_{\lambda,\mu}
t_\lambda t^*_\mu\Big) Q^{\vee F}_\gamma
= Q^{\vee F}_\gamma \Big(\sum_{\substack{\lambda,\mu \in F \\
\gamma(e, d(\lambda)) = \lambda \\
\gamma(e, d(\mu)) = \mu}}
\alpha_{\lambda,\mu} t_\lambda t^*_\mu\Big) Q^{\vee F}_\gamma.
\]
In particular, notice that for $\lambda \in \vee F$,
\begin{equation}\label{eq:Q vee F diagonal mpctn}
Q^{\vee F}_\gamma t_\lambda t^*_\lambda =
\begin{cases}
Q^{\vee F}_\gamma &\text{if $d(\gamma) \ge d(\lambda)$ and $\gamma(e,
d(\lambda)) = \lambda$} \\
0 &\text{otherwise.}
\end{cases}
\end{equation}

We will replace $Q^{\vee F}_\gamma$ with a smaller nonzero projection
$Q_\gamma$ so that the remaining off-diagonal terms are eliminated.
Since $0 < Q_\gamma \le Q^{\vee F}_\gamma$, we will then have
\begin{equation}\label{eq:Q diagonal mpctn}
Q_\gamma t_\lambda t^*_\lambda =
\begin{cases}
Q_\gamma &\text{if $d(\gamma) \ge d(\lambda)$ and $\gamma(e,
d(\lambda)) = \lambda$} \\
0 &\text{otherwise,}
\end{cases}
\end{equation}
which, in conjunction with \eqref{eq:Q vee F diagonal mpctn},
will imply that
\begin{equation}\label{eq:Q norm equality}
\Big\|Q_\gamma \Big(\sum_{\lambda \in F} \alpha_{\lambda,\lambda}
t_\lambda t^*_\lambda \Big)\Big\| =
\Big\lvert \sum_{\substack{\lambda \in F,\, d(\lambda) \leq d(\gamma),
\\ \gamma(e, d(\lambda)) = \lambda}} \alpha_{\lambda,\lambda}
\Big\rvert =
\Big\|Q^{\vee F}_\gamma \Big(\sum_{\lambda \in F}
\alpha_{\lambda,\lambda} t_\lambda t^*_\lambda \Big)\Big\|.
\end{equation}

To produce $Q_\gamma$, we consider pairs $\lambda, \mu \in \vee F$
such that $\gamma(e, d(\lambda)) = \lambda$ and $\gamma(e, d(\mu)) =
\mu$. For each such $(\lambda, \mu)$,
factorise $\gamma$ as $\lambda\lambda' = \gamma=\mu\mu'$, and define
\[
d_{\gamma}(\lambda, \mu) := \big\{\sigma : \sigma =
\delta(d(\lambda^\prime),d(\delta)) \text{ or } \sigma =
\delta(d(\mu^\prime), d(\delta))\text{ for some } \delta \in
\MCE(\lambda^\prime, \mu^\prime)\big\}.
\]
Now $\lambda^\prime$ and $\mu^\prime$ are uniquely determined by
$\lambda$, $\mu$ and $\gamma$, each $\MCE(\lambda^\prime, \mu^\prime)$
is finite, and $\delta(d(\lambda^\prime),d(\delta))$ and
$\delta(d(\mu^\prime),d(\delta))$ are uniquely determined by $\delta
\in \MCE(\lambda^\prime, \mu^\prime)$, so each
$d_\gamma(\lambda,\mu)$ is finite. Let
\[
Q_\gamma := Q^{\vee F}_\gamma
\prod_{\substack{\lambda\not=\mu \in \vee F,\, \gamma(e, d(\lambda)) =
\lambda, \\
\gamma(e, d(\mu)) = \mu,\, \sigma\in d_\gamma(\lambda,\mu)}}
\big(t_\gamma t^*_\gamma - t_{\gamma\sigma} t^*_{\gamma\sigma}\big).
\]
\lemref{lem:Qs nonzero} implies $Q_\gamma > 0$, and $Q_\gamma \le
Q^{\vee F}_\gamma$ by definition, so we have \eqref{eq:Q
diagonal mpctn} and \eqref{eq:Q norm equality}. For
$\lambda,\mu \in \vee F$ with $\lambda\lambda' = \gamma=\mu\mu'$ and
$\lambda \not=
\mu$, we calculate:
\begin{align*}
Q_\gamma t_\lambda t^*_\mu Q_\gamma
&= Q_\gamma (t_{\lambda}(t_{\lambda^\prime}t^*_{\lambda^\prime}
t_{\mu^\prime}
t^*_{\mu^\prime}) t^*_{\mu}) Q_\gamma \\
&= Q_\gamma \Big(t_\lambda \Big(\sum_{\nu \in \MCE(\lambda^\prime,
\mu^\prime)}
t_\nu t^*_\nu \Big) t^*_\mu\Big) Q_\gamma,
\end{align*}
which vanishes because $\nu \in \MCE(\lambda^\prime,
\mu^\prime)$ implies that $\lambda\nu = \gamma\sigma$ for some
$\sigma \in d_\gamma(\lambda,\mu)$.
Thus
\begin{flalign*}
&&\Big\|\sum_{\lambda,\mu \in F} \alpha_{\lambda,\mu} t_\lambda
t^*_\mu \Big\|
&\ge \Big\|Q_\gamma \Big(\sum_{\lambda,\mu \in F}
\alpha_{\lambda,\mu} t_\lambda t^*_\mu \Big) Q_\gamma\Big\| &\\
&&&= \Big\|Q_\gamma \Big(\sum_{\lambda \in F}
\alpha_{\lambda,\lambda} t_\lambda t^*_\lambda \Big) Q_\gamma\Big\|
&\\
&&&= \Big\|Q^{\vee F}_\gamma \Big(\sum_{\lambda \in F}
\alpha_{\lambda,\lambda} t_\lambda t^*_\lambda \Big) \Big\|
\quad\text{by \eqref{eq:Q norm equality}}&\\
&&&= \Big\|\sum_{\lambda \in F} \alpha_{\lambda,\lambda} t_\lambda
t^*_\lambda \Big\|
\qquad\text{by \eqref{eq:norm attained}}. &\qed
\end{flalign*}
\renewcommand{\qed}{}
\end{proof}

\begin{proof}[Proof of \thmref{thm:faithfulness}]
It suffices to show that if $F$ is a finite subset of $E^1$ and
\[
a = \sum_{\lambda,\mu \in F} \alpha_{\lambda,\mu} s_\lambda s_\mu^*,
\]
then $\pi_t(a) = 0$ implies $a=0$. Suppose $\pi_t(a) = 0$. Then
$\pi_t(a^*a) = 0$, $\Phi^B(\pi_t(a^*a))=0$, and \propref{prop:norm
decreasing} implies
that  $\pi_t(\Phi^E(a^*a)) = 0$. Now $\Phi^E(a^*a)$ belongs to $D:=
\lsp\{s_\lambda s_\lambda^* : \lambda \in \vee F\}$, and applying
\propref{prop:partition} to the universal Toeplitz-Cuntz-Krieger $E$-family
$\{s_\lambda\}$ shows that $D$ is a finite-dimensional diagonal
matrix algebra with
matrix units
\[
\Big\{e_{\lambda,\lambda} := s_\lambda
s_\lambda^*\Big(\prod_{\lambda\alpha \in \vee F,\,
d(\alpha) \not= e} (s_{s(\lambda)} -
s_{\lambda\alpha} s_{\lambda\alpha}^*)\Big) :
\lambda \in \vee F\Big\}.
\]
\lemref{lem:Qs nonzero} implies that $\pi_t(e_{\lambda,\lambda})
\not= 0$ for
$\lambda \in \vee F$, so $\pi_t$ is faithful on
$D$.
In particular $\|\Phi^E(a^*a)\| =
\|\pi_t(\Phi^E(a^*a))\| = 0$. \propref{prop:faithful
expectation} now shows that $a^*a = 0$, and hence $a = 0$.
\end{proof}

\section{The $C^*$-algebra of an infinite $k$-graph} \label{sec7}
We show how the finitely-aligned hypothesis, relation (5) of
Definition~\ref{def:new
TCK family}, and the hypothesis
\eqref{eq:faithfulness} in Theorem~\ref{thm:faithfulness} all simplify when the
underlying semigroup is $\NN^k$. We then prove a uniqueness theorem for the
$C^*$-algebras of $k$-graphs in which every vertex receives
infinitely many paths of
every degree.

\subsection{Product systems of graphs over $\NN^k$}

\begin{lemma} \label{lem:NNk f.a.}
Let $(E,\varphi)$ be a product system of graphs over $\NN^k$. Then
$(E,\varphi)$ is finitely aligned if and only if
\begin{equation}\label{faforkgraphs}
\text{$\MCE(\mu, \nu)$ is finite for
every pair $\mu\in E^1_{e_i}$ and $\nu
\in E^1_{e_j}$ with $i\not= j$.}
\end{equation}
\end{lemma}

\begin{proof}
Every finitely aligned system trivially satisfies (\ref{faforkgraphs}). For the
reverse implication, suppose
$E$ satisfies (\ref{faforkgraphs}). Then
$\MCE(\mu,\nu)$ is finite whenever
$\lvert d(\mu) \vee d(\nu)
\rvert \le 2$ . Suppose as an inductive hypothesis that
$\MCE(\mu,\nu)$ is finite whenever $\lvert d(\mu) \vee d(\nu)
\rvert \le n$, and consider $\mu \in E^1_p$, $\nu \in E^1_q$ with $\lvert
p \vee q \rvert = n+1$.

If  the coordinate-wise minimum $p \wedge q$ of $p$ and $q$
is nonzero, then either $\mu(0, p \wedge q) \not= \nu(0, p
\wedge q)$, in which case the factorisation property implies $\MCE(\mu,\nu) =
\emptyset$, or
\[
\MCE(\mu,\nu) = \big\{\mu(0, p \wedge q) \gamma : \gamma \in \MCE(\mu(p
\wedge q, p), \nu(p \wedge q, q))\big\}
\]
is finite by the inductive hypothesis. Thus we may assume that $p
\wedge q = 0$, and hence that $p \vee q = p + q$. If $p \ge q$ or $q \ge p$
then $\MCE(\mu,\nu)$ has at most one element.  So we may further
assume that there exist $i \not= j$ such that $p_i > q_i$
and $q_j > p_j$. Since $p \wedge q = 0$, this implies that $p_j = q_i
= 0$.

Now let $\gamma \in \MCE(\mu,\nu)$. Then $d(\gamma) - e_i = p + q - e_i = (p -
e_i) \vee q$ since $q_i = 0$.
Thus $\gamma_i:= \gamma(0, d(\gamma)-e_i)$ satisfies
\[
\gamma_i(0, p-e_i) = \gamma(0, p-e_i) = \mu(0, p-
e_i) \ \text{ and }\
\gamma_i(0, q) = \gamma(0, q) = \nu,
\]
so $\gamma_i \in
\MCE(\mu(0, p-e_i), \nu)$. Similarly, $\gamma_j:=\gamma(0,
d(\gamma)-e_j)\in\MCE(\mu,
\nu(0, q-e_j))$. But now $p \vee q = d(\gamma_i) + e_i = d(\gamma_j)
+ e_j$, and
since $i
\not= j$, it follows that $d(\gamma) = d(\gamma_i) \vee d(\gamma_j)$.
Furthermore,
$\gamma(0, d(\gamma_i)) = \gamma_i$ and $\gamma(0, d(\gamma_j)) = \gamma_j$, so
$\gamma \in \MCE(\gamma_i, \gamma_j)$. Hence
\[ \lvert\MCE(\mu,\nu)\rvert \le
\sum_{\substack{\gamma_i \in \MCE(\mu(0, p-e_i), \nu) \\ \gamma_j \in
\MCE(\mu,\nu(0, q-e_j))}} \lvert\MCE(\gamma_i, \gamma_j) \rvert.
\]
By the
inductive hypothesis, $\MCE(\mu(0, p-e_i), \nu)$ and $\MCE(\mu,\nu(0,
q-e_j))$
are finite, so the sum has only finitely many terms. Thus
we take $\gamma_i \in \MCE(\mu(0, p-e_i), \nu)$ and $\gamma_j \in
\MCE(\mu,\nu(0, q-e_j))$, and show that $\MCE(\gamma_i, \gamma_j)$ is finite.
If it is nonempty, then the initial segments of degree $(p \vee q) -
e_i - e_j$ of
$\gamma_i$ and
$\gamma_j$ are the same; call it $\beta$, and write
$\gamma_i=\beta\gamma_i'$, $\gamma_j=\beta\gamma_j'$.
Then $d(\gamma_i')=e_i$ and  $d(\gamma_j')=e_j$, so
$|\MCE(\gamma_i, \gamma_j)| =|\MCE(\gamma^\prime_i, \gamma^\prime_j)|$
is finite by (\ref{faforkgraphs}).
\end{proof}

\begin{lemma} \label{lem:N.C. gens do}
Let $(E, \varphi)$ be a finitely aligned product system of graphs over
$\NN^k$. Then a Toeplitz $E$-family $\{t_\lambda\}$ is a Toeplitz-Cuntz-Krieger
$E$-family if and only if
\begin{equation}\label{eq:N.C. gens do}
t^*_\mu t_\nu = \sum_{\mu\alpha = \nu\beta \in \MCE(\mu,\nu)}
t_\alpha t^*_\beta
\end{equation}
for every $\mu \in E^1_{e_i}$ and $\nu \in E^1_{e_j}$ with $s(\mu) =
s(\nu)$ and $i\not= j$.
\end{lemma}

\begin{proof}
Since (\ref{eq:N.C. gens do}) is a special case of Definition~\ref{def:new
TCK family}(5), we have to
show that  (\ref{eq:N.C. gens do}) implies Definition~\ref{def:new
TCK family}(5). If $\lvert
d(\mu)
\vee d(\nu) \rvert \le 2$, this is trivially true. Suppose as an
inductive hypothesis
that
\eqref{eq:finite N.C.} holds whenever $\lvert d(\mu) \vee d(\nu)
\rvert \le n$ for
some $n \ge 2$. Suppose
$\mu \in E^1_p$ and $\nu \in E^1_q$ where $p$ and $q$ satisfy $\lvert
p \vee q \rvert
= n+1$.
We give separate arguments for $p\wedge
q\not=0$ and
$p \wedge q = 0$.

If $p\wedge q\not=0$, then
\begin{align}
t_\mu^* t_\nu
&= t_{\mu(p \wedge q, p)}^* t_{\mu(0, p \wedge q)}^*
t_{\nu(0, p \wedge q)} t_{\nu(p \wedge q, q)} \notag\\
&= \begin{cases}
t_{\mu(p \wedge q, p)}^* t_{\nu(p \wedge q, q)}
&\text{if $\mu(0, p \wedge q) = \nu(0, p \wedge q)$} \label{*}\\
0 &\text{otherwise.}
\end{cases}
\end{align}
The set $\MCE(\mu,\nu)$ is empty unless $\mu(0, p \wedge q) = \nu(0,
p \wedge q)$,
and if so we have
\[
\MCE(\mu,\nu) =
\big\{\mu(0, p \wedge q)\gamma : \gamma \in \MCE(\mu(p \wedge q, p), \nu(p
\wedge q, q))\big\}.
\]
Applying the inductive hypothesis to $(\ref{*})$
gives Definition~\ref{def:new TCK family}(5).

Now suppose $p \wedge q = 0$, or equivalently that $p\vee q=p+q$.
Since $\lvert p \vee
q\rvert\ge 3$, we can assume that $\lvert q
\rvert \ge 2$. If $p \ge q$ then \eqref{eq:finite N.C.} is
trivial, so we may further assume that there exists $i$ such that
$q_i > p_i$, and
then $p \wedge q = 0$ forces $p_i = 0$. In particular,
$\lvert p \vee (q - e_i) \rvert = n$, and the inductive hypothesis gives
\[
t_\mu^* t_\nu =
t_\mu^* t_{\nu(0, q-e_i)} t_{\nu(q-e_i, q)} =
\Big(\sum_{\mu\delta = \nu(0, q-e_i)\varepsilon \in \MCE(\mu,\nu(0, q -
e_i))} t_\delta t_\varepsilon^*\Big) t_{\nu(q-e_i, q)}.
\]
Each $\varepsilon$ appearing in this sum has $d(\varepsilon)=p$, so
$d(\varepsilon)\vee d(\nu(q-e_i,q)=p+e_i$, which has length at most $n$ because
$|q|\geq 2$. Thus we can apply the inductive hypothesis to each summand to get
\begin{equation}\label{ddag}
t_\mu^* t_\nu
= \sum_{\substack{
\mu\delta = \nu(0, q-e_i)\varepsilon \in \MCE(\mu,\nu(0, q - e_i)) \\
\varepsilon\sigma = \nu(q-e_i, q)\tau \in \MCE(\varepsilon, \nu(q- e_i, q))
}}
t_{\delta\sigma}t_\tau^*.
\end{equation}
It remains to show that the pairs $(\delta\sigma,\tau)$ arising in this sum are
precisely the pairs $(\alpha,\beta)$ arising in the right-hand side of
\eqref{eq:finite N.C.}. Given $(\delta\sigma,\tau)$, we certainly have
\[
\mu\delta\sigma = \nu(0, q-e_i)\varepsilon\sigma = \nu(0, q-e_i) \nu(q-e_i,
q) \tau = \nu\tau,
\]
and  $d(\delta\sigma)=d(\delta)+d(\sigma)=q-e_i+e_i=q$, so
$\mu\delta\sigma\in \MCE(\mu,\nu)$. Conversely, given $(\alpha,\beta)$, we take
$\delta:=\alpha(0,q-e_i)$, $\sigma:=\alpha(q-e_i,q)$ and $\tau:=\beta$.
\end{proof}

\begin{lemma}\label{lem:generators do}
Let $E$ be a finitely aligned product system of graphs over $\NN^k$.
Then a Toeplitz
$E$-family $\{t_\lambda\}$ satisfies
\eqref{eq:faithfulness} if and only if
\begin{equation}\label{eq:generators}
\prod^k_{m=1}\Big(t_v - \sum_{\lambda \in G_m} t_\lambda t^*_\lambda
\Big) > 0
\end{equation}
for every choice of finite sets $G_m\subset s^{-1}_{e_m}(v)$.
\end{lemma}

\begin{proof}
The necessity of \eqref{eq:generators} is obvious. Suppose
\eqref{eq:generators} holds, and $R$, $v$ and $F_p$ are
as in \thmref{thm:faithfulness}. For $p\in R$, choose $i_p$ such that
$p_{i_p} > 0$, and for each $m$, set
\[
G_m := \bigcup_{\{p\in R: i_p=m\}}\{\lambda(0, e_m) : \lambda
\in F_p\}.
\]
Then each $G_m$ is a finite subset of $s^{-1}_{e_m}(v)$, and
\begin{align*}
\prod_{p\in R}\Big(t_v - \sum_{\lambda\in F_p}
t_\lambda t^*_\lambda\Big)
&\ge \prod_{p\in R}\Big(t_v -
\sum_{\lambda\in F_p} t_{\lambda(0, e_{i_p})} t^*_{\lambda(0,e_{i_p})}
\Big) \\
&= \prod^k_{m = 1}\Big(t_v - \sum_{\mu \in G_m} t_\mu t^*_\mu \Big),
\end{align*}
which is nonzero by \eqref{eq:generators}.
\end{proof}

\subsection{The $C^*$-algebra of an infinite $k$-graph}
If $(\Lambda, d)$ is a $k$-graph, and $\lambda, \mu \in \Lambda$, we
regard $\MCE(\lambda,\mu) \subset (E_\Lambda)^1$ as a subset of
$\Lambda$. In view of Lemma~\ref{lem:N.C. gens do},  we say that $\Lambda$
is finitely aligned if $\MCE(\lambda,\mu)$ is finite whenever
$d(\lambda)=e_i$ and
$d(\mu)={e_j}$. By a
Toeplitz-Cuntz-Krieger $\Lambda$-family we mean a Toeplitz-Cuntz-Krieger
$E_\Lambda$-family. If $\Lambda$ has no sources, so that the graphs in
$E_\Lambda$ have no sinks, then we define a Cuntz-Krieger
$\Lambda$-family to be a
Cuntz-Pimsner
$E_\Lambda$-family. We have only made this last definition for
$k$-graphs without
sources to avoid clashing with the definitions given for row-finite graphs in
\cite{RSY}; for row-finite $k$-graphs without sources, therefore, our
$C^*(E_\Lambda)$
coincides with the graph algebra $C^*(\Lambda)$ used in \cite{KP} and
\cite{RSY}.

Recall that $\Lambda^n(v):=\{\lambda\in\Lambda: d(\lambda)=n\text{ and
}\cod(\lambda)=v\}$. If
$\lvert\Lambda^{e_i}(v)\rvert =
\infty$ for every
$v
\in
\Lambda^0$, and every $1 \le i \le k$, then conditions (6)  and (4) of
Definition~\ref{def:new TCK family} are equivalent, so
Theorem~\ref{thm:faithfulness} gives a uniqueness theorem for
$C^*(\Lambda)$.

\begin{cor}
Let $(\Lambda,d)$ be a finitely aligned $k$-graph such that
$\lvert\Lambda^{e_i}(v)\rvert = \infty$ for every $v \in \Lambda^0$
and $1 \le i \le k$. Let $\{t_\lambda : \lambda \in \Lambda\}$ be a
Cuntz-Krieger $\Lambda$-family such that $t_v \not= 0$ for all $v \in
\Lambda^0$. Then the representation $\pi_t$ of
$C^*(\Lambda):=C^*(E_\Lambda)$ is
faithful.
\end{cor}

\begin{proof}
That each $\lvert\Lambda^{e_i}(v)\rvert = \infty$ implies both that
$C^*(E_\Lambda)=\Tt C^*(E_\Lambda)$, and that $\Lambda$ has no sources, so that
$C^*(\Lambda):=C^*(E_\Lambda)$. \lemref{lem:NNk f.a.}
implies that $(E_\Lambda,
\varphi_\Lambda)$ is finitely aligned.
To establish
\eqref{eq:faithfulness}, we fix $v \in \Lambda^0$ and finite sets $G_m
\subset \Lambda^{e_m}(v)$ for $1 \le m \le k$. By \lemref{lem:generators
do}, it suffices to show that
\[
\prod^k_{m=1}\Big(t_v - \sum_{\lambda \in G_m} t_\lambda t^*_\lambda
\Big) > 0.
\]

We shall construct paths $\mu_m\in
\Lambda(v)$ of degree $\sum_{i=1}^m e_i$ for $m\leq k$ such that
$\mu_m(0,e_i)$ does
not belong to
$G_i$ for
$1\leq i\leq m$. We take $\mu_1$ to be any edge of degree $e_1$ which is not in
$G_1$. If we have
$\mu_m$, then because  the set $\Lambda^{e_{m+1}}(\dom(\mu_m))$ is
infinite, there is
a path $\mu_{m+1}=\mu_m\alpha$ of degree $\sum_{i=1}^{m+1} e_i$
which is not in the
finite set
$\bigcup_{\lambda
\in G_{m+1}}\MCE(\mu_m,
\lambda)$. Then $\mu_{m+1}(0,e_i)=\mu_m(0,e_i)$ is not in $G_i$ for
$i\leq m$, and
$\mu_{m+1}(0,e_{m+1})$ cannot be in $G_{m+1}$ because $\mu_{m+1}\in
\MCE(\mu_m,\mu(0,e_{m+1}))$.

Now for $\lambda \in G_i$, we have $\MCE(\lambda,\mu_k) = \emptyset$,
and relation
(5) of Definition~\ref{def:new TCK family} in the form (\ref{altNicacov}) gives
$t_\lambda t^*_\lambda t_{\mu_k} t^*_{\mu_k} =0$.
  Thus
\[
\prod^k_{m=1}\Big(t_v - \sum_{\lambda \in G_m} t_\lambda
t^*_\lambda \Big) t_{\mu_k} t^*_{\mu_k}
= t_{\mu_k} t^*_{\mu_k},
\]
which is nonzero because $t_{\mu_k}^*t_{\mu_k}=t_{s(\mu_k)}$ is nonzero. Since
$\ZZ^k$ is amenable, the result now follows from
Theorem~\ref{thm:faithfulness}.
\end{proof}


\begin{thebibliography}{20}


\bibitem{C} J. Cuntz, {\em Simple $C^*$-algebras generated by
isometries},  Comm.
Math. Phys. {\bf 57} (1977), 173--185.


\bibitem{CK1} J. Cuntz and W. Krieger, {\em A class of $C^*$-algebras
and topological Markov chains}, Invent. Math. {\bf 56} (1980),
251--268.

\bibitem{D} D. Drinen, {\em Viewing AF-algebras as graph algebras},
Proc. Amer. Math.
Soc. {\bf 128} (2000),  1991--2000.

\bibitem{F} N. J. Fowler, {\em Compactly aligned product systems, and
generalizations
of
$\Oo_\infty$}, Inter. J. Math. {\bf 10} (1999), 721--738.

\bibitem{F99} N. J. Fowler, {\em Discrete product systems of Hilbert
bimodules}, Pacific J. Math, to appear.

\bibitem{FLR} N. J. Fowler, M. Laca and I. Raeburn, {\em The
$C^*$-algebras of infinite graphs}, Proc. Amer. Math. Soc. {\bf 128}
(2000), 2319--2327.

\bibitem{FR} N. J. Fowler and I. Raeburn, {\em The Toeplitz algebra of
a Hilbert bimodule}, Indiana Univ. Math. J. {\bf 48} (1999), 155--181.

\bibitem{FS} N. J. Fowler and A. Sims, {\em Product systems over
right-angled Artin semigroups}, Trans. Amer. Math. Soc. {\bf 354}
(2002), 1487--1509.

\bibitem{HS} J. H. Hong and W. Szyma{\' n}ski, {\em Quantum spheres
and projective
spaces as graph algebras}, preprint, U. of Newcastle, January 2001.

\bibitem{KPW} T. Kajiwara, C. Pinzari and Y. Watatani, {\em Ideal
structure and simplicity of the $C^*$-algebras generated by Hilbert
bimodules}, J. Funct. Anal. {\bf 159} (1998), 295--322.

\bibitem{KP} A. Kumjian and D. Pask, {\em Higher rank graph
$C^*$-algebras}, New York J. Math {\bf 6} (2000), 1--20.

\bibitem{KPR} A. Kumjian, D. Pask and I. Raeburn, {\em Cuntz-Krieger
algebras of directed graphs,} Pacific J. Math. {\bf 184} (1998),
161--174.

\bibitem{KPRR} A. Kumjian, D. Pask, I. Raeburn and J. Renault, {\em
Graphs, groupoids, and Cuntz-Krieger algebras,} J. Funct. Anal. {\bf
144} (1997), 505--541.

\bibitem{LR} M. Laca and I. Raeburn, {\em Semigroup crossed products
and the Toeplitz
algebras of nonabelian groups}, J. Funct. Anal. {\bf 139} (1996), 415--440.

\bibitem{N92} A. Nica, {\em $C^*$-algebras generated by isometries and
Wiener-Hopf operators,} J. Operator Theory {\bf 27} (1992), 17--52.

\bibitem{P97} M. V. Pimsner, {\em A class of $C^*$-algebras
generalizing both Cuntz-Krieger algebras and crossed products by
$\ZZ$}, Fields Institute Comm. {\bf 12} (1997), 189--212.

\bibitem{RSY} I. Raeburn, A. Sims and T. Yeend, {\em Higher-rank
graphs and their $C^*$-algebras,} Proc. Edinburgh Math. Soc., to appear.

\bibitem{RS} G. Robertson and T. Steger, {\em Affine buildings, tiling
systems and higher rank Cuntz-Krieger algebras,} J. reine angew. Math.
{\bf 513} (1999), 115--144.

\bibitem{Sz} W. Szyma{\' n}ski, {\em The range of $K$-invariants for
$C^*$-algebras
of infinite graphs}, Indiana Univ. Math. J., to appear.

\end{thebibliography}
\end{document}